\documentclass[12pt]{article}
\usepackage{amssymb,amsmath,amsopn,amsfonts}
\usepackage[dvips]{color}
\usepackage{colordvi,multicol}
\usepackage{epsf}
\newfam\msbfam
\font\tenmsb=msbm10 \textfont\msbfam=\tenmsb \font\sevenmsb=msbm7
\scriptfont\msbfam=\sevenmsb \font\fivemsb=msbm5
\scriptscriptfont\msbfam=\fivemsb

\def\th#1{\vspace{1mm}\noindent{\bf #1}\quad}

\def\proof{\vspace{1mm}\noindent{\it Proof}\quad}

\numberwithin{equation}{section}

\def\bc{\begin{center}}
\def\ec{\end{center}}
\def\no{\noindent}
\def\hang{\hangindent\parindent}
\def\textindent#1{\indent\llap{\qquad #1\ \ \enspace}\ignorespaces}
\def\ref{\par\hang\textindent}

\begin{document}
\bigbreak

\title{ {\bf Restricted Markov uniqueness for the stochastic quantization of $P(\Phi)_2$ and its applications
\thanks{Research supported in part  by NSFC (No.11301026, No.11401019), Key Lab of Random Complex Structures and
Data Science, Chinese Academy of Sciences (Grant No. 2008DP173182) and DFG through IRTG 1132 and CRC 701}\\} }
\author{{\bf Michael R\"{o}ckner}$^{\mbox{c},}$, {\bf Rongchan Zhu}$^{\mbox{a},}$\thanks{
 Corresponding author} {\bf Xiangchan Zhu}$^{\mbox{b},}$,
\date {}
\thanks{E-mail address:  roeckner@math.uni-bielefeld.de(M. R\"{o}ckner), zhurongchan@126.com(R. C. Zhu), zhuxiangchan@126.com(X. C. Zhu)}\\ \\
$^{\mbox{a}}$Department of Mathematics, Beijing Institute of Technology, Beijing 100081,
 China,\\
$^{\mbox{b}}$School of Science, Beijing Jiaotong University, Beijing 100044, China\\
$^{\mbox{c}}$ Department of Mathematics, University of Bielefeld, D-33615 Bielefeld, Germany,}

\maketitle

\noindent {\bf Abstract}

In this paper we obtain restricted Markov uniqueness of the generator and uniqueness of probabilistically weak solutions for the stochastic quantization problem in both the finite and  infinite volume case by clarifying the precise relation between the solutions to the stochastic quantization problem  obtained by the Dirichlet form approach and those obtained in [DD03] and in [MW15].  We  prove that the solution $X-Z$, where $X$ is obtained by the Dirichlet form approach in [AR91] and $Z$ is the corresponding O-U process,  satisfies the corresponding  shifted equation (see (1.4) below). Moreover, we obtain that the infinite volume $P(\Phi)_2$ quantum  field is an invariant measure for the $X_0=Y+Z$, where $Y$ is the unique solution to the shifted equation.

\vspace{1mm}
\no{\footnotesize{\bf 2000 Mathematics Subject Classification AMS}:\hspace{2mm} 60H15, 82C28}%Five letters(Mathematics Subject Classification 2000)
 \vspace{2mm}

\no{\footnotesize{\bf Keywords}:\hspace{2mm}  stochastic quantization problem, Dirichlet forms, space-time white noise, Wick power}% At least three keywords

\section{Introduction}

In this paper we analyze stochastic quantization equations on $\mathbb{T}^2$ and on $\mathbb{R}^2$:
So, let $H=L^2(\mathbb{T}^2)$ or $L^2(\mathbb{R}^2)$ and consider
\begin{equation}\aligned dX=&(AX-:p(X):)dt+dW(t),\\X(0)=&z,\endaligned\end{equation}
where $A:D(A)\subset H\rightarrow H$ is the linear operator
$$A\phi=\Delta \phi-\phi,\quad p(\phi)=\sum_{n=1}^{2N}na_n\phi^{n-1},$$
 where $a_{2N}>0$ and $:p(\phi):$ means the renormalization of $p(\phi)$ whose definition we will give in Section 3 and Section 4. $W$ is a cylindrical $\mathcal{F}_t$-Wiener process defined on a probability space $(\Omega,\mathcal{F},P)$ with a normal filtration $(\mathcal{F}_t)_{t\geq0}$.

This equation arises in stochastic quantization of  Euclidean quantum field
theory. Heuristically, (1.1) has an invariant measure $\nu$ defined as $$\nu(d\phi)=ce^{-\int:q(\phi):dx}\mu(d\phi),$$ where $q(\phi)=\sum_{n=0}^{2N}a_n\phi^n$, $c$ is a normalization constant and $\mu$ is the Gaussian free field. $\nu$ is called the $P(\Phi)_2$-quantum field.  There have been many approaches to the problem of
giving a meaning to the above heuristic measure for the two dimensional case and the three dimensional case (see  [GRS75], [GlJ86] and references
therein).
In [PW81]  Parisi and Wu proposed a program for Euclidean quantum field
theory of getting Gibbs states of classical
statistical mechanics as limiting distributions of stochastic processes,  especially as solutions to non-linear stochastic differential
equations. Then one can use the stochastic differential equations to study the properties of the Gibbs states. This
procedure is called stochastic field quantization (see [JLM85]). The $P(\Phi)_2$ model is the simplest non-trivial Euclidean quantum field (see [GlJ86] and the reference therein). The issue of the stochastic quantization of the $P(\Phi)_2$ model is to solve the equation (1.1).

The main difficulty in this case is that $W$
and hence the solutions are so singular that the non-linear term is not well-defined in the classical sense.
 In [AR91]  weak solutions to (1.1) have been constructed by using the Dirichlet form approach in the finite  and infinite volume case.
 However, \emph{Markov uniqueness} for the corresponding generator $(L, D)$ has been an open problem for many years. Here $D$ is the "minimal" domain contained in the domain of the generator. Consider a measure $\nu$ on a Banach space $E$. The problem of Markov uniqueness is whether
 there exists exactly one negative definite self-adjoint operator $L^\nu$ on $L^2(E;\nu)$ which extends $(L, D)$ and is a Dirichlet operator, i.e. $T_t:=e^{tL^\nu}$ is sub-Markovian. The latter property is equivalent to the quadratic form given by $L^\nu$ on $L^2(E;\nu)$ being a\emph{ Dirichlet form}.
 Then Markov uniqueness is equivalent to the fact that there exists exactly one  Dirichlet form whose generator extends  $(L, D)$.
  This problem is
completely solved in the finite dimensional case in [RZ94] where Markov uniqueness
was obtained under the most general conditions. The situation is quite different in the infinite dimensional case. We refer to [ARZ93a], [ARZ93b], [LR98], [KR07], [AKR12] for the best results in this direction
known so far. In these papers Markov uniqueness has been obtained for a modified  stochastic quantization equation
$$\aligned dX=&(-\Delta+1)^{-\varepsilon}(AX-:p(X):)dt+(-\Delta+1)^{-\frac{\varepsilon}{2}}dW(t),\endaligned$$
with $\varepsilon>0$. However, Markov uniqueness for the case that $\varepsilon=0$ is still an open problem.

In this paper we study Markov uniqueness for the  operator associated with the stochastic quantization problem in both the finite volume case and the infinite volume case and obtain the restricted Markov uniqueness of the  operator, i.e.  there exists exactly one \emph{quasi-regular Dirichlet form} whose generator extends $(L, D)$ (see Theorem 3.12, Theorem 4.10).

 This problem is also related to the uniqueness of the  martingale problem for $(L, D)$, i.e. whether
 there exists exactly one (up to $\nu$-equivalence defined in Section 3) strong Markov process solving the martingale problem for $(L, D)$, and the uniqueness of probabilistically weak solution to (1.1). In this paper we also obtain that there exists exactly one  (up to  $\nu$-equivalence defined in Section 3) probabilistically weak solution to (1.1) (see Theorem 3.12 and Theorem 4.10).

 We obtain Markov uniqueness in the restricted sense and  the uniqueness of the  probabilistically weak solution to (1.1) by studying the relations between the solutions to the stochastic quantization problem obtained by the Dirichlet form approach and those obtained in [DD03] and in [MW15]. In fact,  (1.1) has been studied by many authors:
 In [MR99] the stationary solution to (1.1) has also been considered in their general theory of martinglae solutions for stochastic partial differential equaitons; In [DD03] Da Prato and Debussche define the Wick powers of solutions to the stochastic heat equation in the paths space and study a shifted equation instead of (1.1) in the finite volume case. They split the unknown $X$ into two parts:
$X=Y_1+Z_1$, where $Z_1(t)=\int_{-\infty}^te^{(t-s)A}dW(s)$. Observe that $Y_1$ is much smoother than $X$ and that in the stationary case
\begin{equation}:X^k:=\sum_{l=0}^kC_k^lY_1^l:Z_1^{k-l}:,\end{equation}
with $C_k^l=\frac{k!}{l!(k-l)!}$,
which motivate them  to consider the following shifted equation:
\begin{equation}\aligned \frac{dY_1}{dt}=&AY_1-\sum_{k=1}^{2N}ka_k\sum_{l=0}^{k-1} C_{k-1}^lY_1(t)^l:Z_1(t)^{k-1-l}:\\Y_1(0)=&z-Z_1(0),\endaligned\end{equation}
and obtain  local existence and uniqueness of the solution $Y_1$ to (1.3) by a fixed point argument. By using the invariant measure $\nu$ they obtain a global solution to (1.1) by defining $X=Y_1+Z$ starting from almost every starting point. In [MW15] the authors consider the following equation with $N=2$ instead of (1.3):
 \begin{equation}\aligned \frac{dY}{dt}=&AY-\sum_{k=1}^{2N}ka_k\sum_{l=0}^{k-1} C_{k-1}^lY(t)^l:\bar{Z}(t)^{k-1-l}:\\Y(0)=&0,\endaligned\end{equation}
 where $\bar{Z}(t)=e^{tA}z+\int_0^te^{(t-s)A}dW(s)$. We call (1.4) \emph{the shifted equation} for short.
 They obtain global existence and uniqueness of the solution to (1.4) directly from every starting point both in the finite and  infinite volume case.
Actually, (1.3) is equivalent to (1.4). For the solution $Y_1$ to (1.3), defining $Y(t):=Y_1(t)+e^{tA}Z_1(0)-e^{tA}z$, we can easily check that $Y$ is a solution to (1.4) by using the binomial formula (3.1) below.

It is natural to ask whether the unique solution obtained by the methods in [DD03] and [MW15]  satisfies the original equation (1.1) and has $\nu$ as  an invariant measure. Furthermore, it is a priori far from being clear what is the relation between the solutions obtained by the Dirichlet form approach and the solution obtained in [DD03] and in [MW15].  In this paper we study this problem and we prove that  $X-\bar{Z}$, where $X$ is obtained by the Dirichlet form approach in [AR91] and $\bar{Z}(t)=\int_0^te^{(t-s)A}dW(s)+e^{tA}z$, also satisfies the shifted equation (1.4). We emphasize that it is not obvious that $X-\bar{Z}$ satisfies the shifted equation (1.4) since  (1.2) only holds in the stationary case and we do not know whether the marginal distribution of the solution is absolutely continuous with respect to $\nu$. However, by using Dirichlet form theory we can solve this problem and obtain the desired results (see Theorem 3.9 and Theorem 4.8).

 Moreover, we obtain that the $P(\Phi)_2$ quantum  field $\nu$ is an invariant measure for the process $X_0=Y+\bar{Z}$, where $Y$ is the unique solution to the shifted equation (1.4). As a consequence, we deduce uniqueness of probabilistically weak solutions to (1.1) and Markov uniqueness for the corresponding generator in the restricted sense in both the finite  and infinite volume case (see Theorem 3.12 and Theorem 4.10). We also emphasize that the $P(\Phi)_2$ field is not absolutely continuous with respect to Gaussian measure in the infinite volume case. This makes it more difficulty to analyze the support of $\nu$. Here we use [GlJ86] and techniques from Dirichlet form theory to  solve this problem (see Theorem 4.7).

We also want to mention that recently there
has arisen a renewed interest in SPDE’s related to such problems, particularly in
connection with Hairer's theory of regularity structures [Hai14] and related
work by Imkeller, Gubinelli, Perkowski in [GIP13]. By using these theories one can obtain local existence and uniqueness of solution to (1.1) in the three dimensional case (see [Hai14, CC13]). In a forthcoming paper we also prove ergodicity for the finite volume
case. To the best
of our knowledge, this is still an open problem in the periodic case
(i.e. on the torus).
In the infinte volume case this has been studied in [AKR97].

This paper is organized as follows: In Section 2 we collect some results related to Besov  and weighted Besov spaces. In Section 3 we consider the finite volume case and prove that  the solution obtained by Dirichlet form theory  satisfies the shifted equation. Moreover, we obtain Markov uniqueness in the restricted sense and uniqueness of the  probabilistically weak solutions to (1.1).  In Section 4 we prove that all the results also hold in the infinite volume case.

\section{Preliminary}
In the following we recall the definitions  of Besov spaces. For a general introduction to the theory we refer to [BCD11, Tri78, Tri06].
 The space of real valued infinitely differentiable functions of compact support is denoted by $\mathcal{D}(\mathbb{R}^d)$ or $\mathcal{D}$. The space of Schwartz functions is denoted by $\mathcal{S}(\mathbb{R}^d)$. Its dual, the space of tempered distributions, is denoted by $\mathcal{S}'(\mathbb{R}^d)$. The Fourier transform and the inverse Fourier transform are denoted by $\mathcal{F}$ and $\mathcal{F}^{-1}$, respectively.

 Let $\chi,\theta\in \mathcal{D}$ be nonnegative radial functions on $\mathbb{R}^d$, such that

i. the support of $\chi$ is contained in a ball and the support of $\theta$ is contained in an annulus;

ii. $\chi(z)+\sum_{j\geq0}\theta(2^{-j}z)=1$ for all $z\in \mathbb{R}^d$.

iii. $\textrm{supp}(\chi)\cap \textrm{supp}(\theta(2^{-j}\cdot))=\emptyset$ for $j\geq1$ and $\textrm{supp}\theta(2^{-i}\cdot)\cap \textrm{supp}\theta(2^{-j}\cdot)=\emptyset$ for $|i-j|>1$.

We call such $(\chi,\theta)$ dyadic partition of unity, and for the existence of dyadic partitions of unity see [BCD11, Proposition 2.10]. The Littlewood-Paley blocks are now defined as
$$\Delta_{-1}u=\mathcal{F}^{-1}(\chi\mathcal{F}u)\quad \Delta_{j}u=\mathcal{F}^{-1}(\theta(2^{-j}\cdot)\mathcal{F}u).$$

For $\alpha\in\mathbb{R}$, $p,q\in [1,\infty]$, $u\in\mathcal{D}$ we define
$$\|u\|_{B^\alpha_{p,q}}:=(\sum_{j\geq-1}(2^{j\alpha}\|\Delta_ju\|_{L^p})^q)^{1/q},$$
with the usual interpretation as $l^\infty$ norm in case $q=\infty$. The Besov space $B^\alpha_{p,q}$ consists of the completion of $\mathcal{D}$ with respect to this norm and the H\"{o}lder-Besov space $\mathcal{C}^\alpha$ is given by $\mathcal{C}^\alpha(\mathbb{R}^d)=B^\alpha_{\infty,\infty}(\mathbb{R}^d)$. For $p,q\in [1,\infty)$,
$$B^\alpha_{p,q}(\mathbb{R}^d)=\{u\in\mathcal{S}'(\mathbb{R}^d):\|u\|_{B^\alpha_{p,q}}<\infty\}.$$
$$\mathcal{C}^\alpha(\mathbb{R}^d)\varsubsetneq \{u\in\mathcal{S}'(\mathbb{R}^d):\|u\|_{\mathcal{C}^\alpha(\mathbb{R}^d)}<\infty\}.$$
We point out that everything above and everything that follows can be applied to distributions on the torus (see [S85, SW71]). More precisely, let $\mathcal{S}'(\mathbb{T}^d)$ be the space of distributions on $\mathbb{T}^d$.  Besov spaces on the torus with general indices $p,q\in[1,\infty]$ are defined as
the completion of $\mathcal{D}$ with respect to the norm $$\|u\|_{B^\alpha_{p,q}(\mathbb{T}^d)}:=(\sum_{j\geq-1}(2^{j\alpha}\|\Delta_ju\|_{L^p(\mathbb{T}^d)})^q)^{1/q},$$
and the H\"{o}lder-Besov space $\mathcal{C}^\alpha$ is given by $\mathcal{C}^\alpha=B^\alpha_{\infty,\infty}(\mathbb{T}^d)$.  We write $\|\cdot\|_{\alpha}$ instead of $\|\cdot\|_{B^\alpha_{\infty,\infty}(\mathbb{T}^d)}$ in the following for simplicity.  For $p,q\in[1,\infty)$
$$B^\alpha_{p,q}(\mathbb{T}^d)=\{u\in\mathcal{S}'(\mathbb{T}^d):\|u\|_{B^\alpha_{p,q}(\mathbb{T}^d)}<\infty\}.$$
\begin{equation}\mathcal{C}^\alpha\varsubsetneq \{u\in\mathcal{S}'(\mathbb{T}^d):\|u\|_{\alpha}<\infty\}.\end{equation}

Here we choose Besov spaces as  completions of smooth functions with compact support, which ensures that the Besov spaces are separable which has a lot of advantages for our analysis below.
 \vskip.10in

{\textbf{Weighted Besov spaces}}

In the following we recall the definitions and some properties of weighted Besov spaces, which are used for analyzing the regularity of the distributions in the infinite volume case.  For a general introduction to these theories we refer to [Tri06].

For $\sigma\in\mathbb{R}$ we let $w(x)=(1+|x|^2)^{-\sigma/2}$. For $\alpha\in \mathbb{R}$,  $p,q\in [1,\infty]$, we define the weighted Besov norm for $u\in \mathcal{D}$,
$$\|u\|_{\hat{\mathcal{B}}^{\alpha,\sigma}_{p,q}}
:=(\sum_{j\geq-1}(2^{j\alpha}\|\Delta_ju \|_{L^p(wdx)})^q)^{1/q},$$
with the usual interpretation as $l^\infty$ norm in case $q=\infty$ and $\|f\|_{L^\infty(wdx)}=\|wf\|_{L^\infty(dx)}$. The Besov space $\hat{\mathcal{B}}^{\alpha,\sigma}_{p,q}(\mathbb{R}^d)$ consists of the completion of $\mathcal{D}$ with respect to this norm. For $p, q\in[1,\infty)$,
$$\hat{\mathcal{B}}^{\alpha,\sigma}_{p,q}(\mathbb{R}^d)=\{u\in\mathcal{S}'(\mathbb{R}^d):\|u\|_{\hat{\mathcal{B}}^{\alpha,\sigma}_{p,q}}
<\infty\}.$$
By [Tri83, Theorem 9.2.1] we can view functions on the torus $\mathbb{T}^d$ as periodic functions on $\mathbb{R}^d$  and have that for $\alpha\in\mathbb{R},\sigma>0$, $f\in \mathcal{C}^\alpha$
\begin{equation}\|f\|_{\alpha}\backsimeq\|f\|_{\hat{\mathcal{B}}^{\alpha,\sigma}_{\infty,\infty}}.\end{equation}

\textbf{Wavelet analysis}

We will also use wavelet analysis to determine the regularity of a distribution in a Besov space.
In the following we briefly summarize wavelet analysis below; we refer to  work of Meyer [Mey92],
Daubechies [Dau88] and [Tri06] for more details on wavelet analysis. For every $r > 0$, there exists a compactly supported function
$\varphi\in C^r(\mathbb{R})$ such that:

1. We have $\langle \varphi(\cdot),\varphi(\cdot-k)\rangle=\delta_{k,0}$ for every $ k\in \mathbb{Z}$,

2. There exist ˜$\tilde{a}_k,k\in\mathbb{Z}$ with only finitely many non-zero values, and such that
$\varphi(x)=\sum_{k\in\mathbb{Z}}\tilde{a}_k\varphi(2x-k)$ for every $x\in\mathbb{R}$,

3. For every polynomial $P$ of degree at most $r$ and for every $x\in\mathbb{R}$,
$\sum_{k\in\mathbb{Z}}\int P(y)\varphi(y-k)dy \varphi(x-k) = P(x)$.

Given such a function $\varphi$, we define for every $ x\in \mathbb{R}^d$ the recentered and
rescaled function $\varphi^n_x$ as follows
$$\varphi^n_x(y) := \Pi_{i=1}^d2^{\frac{n}{2}}\varphi(2^n(y_i-x_i)).$$
Observe that this rescaling preserves the $L^2$-norm. We let $V_n$ be the subspace of
$L^2(\mathbb{R}^d)$ generated by $\{\varphi^n_x: x\in \Lambda_n\}$, where
$$\Lambda_n := \{( 2^{−n}k_1,...,2^{−n}k_d) : k_i\in\mathbb{Z}\}.$$
An important property of wavelets is the existence of a finite set $\Psi$	 of compactly
supported functions in $C^r$ such that, for every $n\geq0$, the orthogonal complement
of $V_n$ inside $V_{n+1}$ is given by the linear span of all the $\psi^n_x, x\in \Lambda_n, \psi\in \Psi$.
 For every $n\geq0$
$$\{\varphi^n_x, x\in\Lambda_n\}\cup \{\psi^m_x: m\geq n, \psi\in\Psi, x\in \Lambda_m\},$$
forms an orthonormal basis of $L^2(\mathbb{R}^d)$.
This wavelet analysis allows one to identify a countable collection of conditions
that determine the regularity of a distribution.

 Setting ${\Psi}_\star={\Psi}\cup\{\varphi\}$, by [Tri06, Theorem 6.15] we know that for $p\in(1,\infty)$, $\alpha\in\mathbb{R}$, $f\in \hat{\mathcal{B}}^{\alpha,\sigma}_{p,p}$
\begin{equation}\|f\|^p_{\hat{\mathcal{B}}^{\alpha,\sigma}_{p,p}}\leq C\sum_{n=0}^\infty 2^{n(\alpha-d/p+1)p}\sum_{\psi\in\Psi_\star}\sum_{x\in\Lambda_n}|\langle f,\psi^n_{x}\rangle|^pw(x),\end{equation}
and \begin{equation}\|f\|^p_{\hat{\mathcal{B}}^{\alpha,\sigma}_{\infty,\infty}}\leq C\sum_{n=0}^\infty 2^{n(\alpha+1)p}\sum_{\psi\in\Psi_\star}\sum_{x\in\Lambda_n}|\langle f,\psi^n_{x}\rangle|^pw(x)^p.\end{equation}
 \vskip.10in
\textbf{Estimates on the torus}

 Set $\Lambda= (-A)^{\frac{1}{2}}$. For $s\geq0, p\in [1,+\infty]$ we use $H^{s}_p$ to denote the subspace of $L^p(\mathbb{T}^d)$, consisting of all  $f$   which can be written in the form $f=\Lambda^{-s}g, g\in L^p(\mathbb{T}^d)$ and the $H^{s}_p$ norm of $f$ is defined to be the $L^p$ norm of $g$, i.e. $\|f\|_{H^{s}_p}:=\|\Lambda^s f\|_{L^p(\mathbb{T}^d)}$.

\vskip.10in
To study (1.1) in the finite volume case, we  will need several important properties of Besov spaces on the torus and we recall the following Besov embedding theorems on the torus first (c.f. [Tri78, Theorem 4.6.1], [GIP13, Lemma 41]):
\vskip.10in
 \th{Lemma 2.1} (i) Let $1\leq p_1\leq p_2\leq\infty$ and $1\leq q_1\leq q_2\leq\infty$, and let $\alpha\in\mathbb{R}$. Then $B^\alpha_{p_1,q_1}(\mathbb{T}^d)$ is continuously embedded in $B^{\alpha-d(1/p_1-1/p_2)}_{p_2,q_2}(\mathbb{T}^d)$.

 (ii) Let $s\geq0$, $1<p<\infty$, $\epsilon>0$. Then
 $ H^{s+\epsilon}_p\subset B^{s}_{p,1}(\mathbb{T}^d)\subset B^{s}_{1,1}(\mathbb{T}^d)$.

 (iii) Let $1\leq p_1\leq p_2<\infty$ and let $\alpha\in\mathbb{R}$. Then $H^\alpha_{p_1}$ is continuously embedded in $H^{\alpha-d(1/p_1-1/p_2)}_{p_2}$.

Here  $\subset$ means that the embedding is continuous and dense.

\vskip.10in

We recall the following Schauder estimates, i.e. the smoothing effect of the heat flow, for later use.

\vskip.10in
\th{Lemma 2.2}([GIP13, Lemma 47]) (i) Let $u\in B^{\alpha}_{p,q}(\mathbb{T}^d)$ for some $\alpha\in \mathbb{R}, p,q\in [1,\infty]$. Then for every $\delta\geq0$
$$\|e^{tA}u\|_{B^{\alpha+\delta}_{p,q}(\mathbb{T}^d)}\lesssim t^{-\delta/2}\|u\|_{B^{\alpha}_{p,q}(\mathbb{T}^d)}.$$
(ii) Let $\alpha\leq \beta\in\mathbb{R}$. Then
$$\|(1-e^{tA})u\|_{\alpha}\lesssim t^{\frac{\beta-\alpha}{2}}\|u\|_{\beta}.$$

\vskip.10in
One can  extend the multiplication on suitable Besov spaces and also  have the duality properties of Besov spaces from [Tri78, Chapter 4]:
\vskip.10in

\th{Lemma 2.3} (i) The bilinear map $(u; v)\mapsto uv$
extends to a continuous map from $\mathcal{C}^\alpha\times \mathcal{C}^\beta$ to $\mathcal{C}^{\alpha\wedge\beta}$ if and only if $\alpha+\beta>0$.

(ii) Let $\alpha\in (0,1)$, $p,q\in[1,\infty]$, $p'$ and $q'$ be their conjugate exponents, respectively. Then the mapping  $(u; v)\mapsto \int uvdx$  extends to a continuous bilinear form on $B^\alpha_{p,q}(\mathbb{T}^d)\times B^{-\alpha}_{p',q'}(\mathbb{T}^d)$.

\vskip.10in
We recall the following interpolation inequality and  multiplicative inequality for the elements in $H^s_p$, which is required for the a-priori estimate in the proof of Theorem 3.10: (cf. [Tri78, Theorem 4.3.1], [Re95, Lemma A.4], [RZZ15a, Lemma 2.1]):
 \vskip.10in

\th{Lemma 2.4} (i)  Suppose that $s\in (0,1)$ and $p\in (1,\infty)$. Then for $u\in H^1_p$
$$\|u\|_{H^s_p}\lesssim \|u\|_{L^p(\mathbb{T}^d)}^{1-s}\|u\|_{H^1_p}^s.$$

(ii) Suppose that $s>0$ and $p\in (1,\infty)$. If $u,v\in C^{\infty}(\mathbb{T}^2)$ then
$$\|\Lambda^s(uv)\|_{L^p(\mathbb{T}^d)}\leq C(\|u\|_{L^{p_1}(\mathbb{T}^d)}\|\Lambda^sv\|_{L^{p_2}(\mathbb{T}^d)}+\|v\|_{L^{p_3}(\mathbb{T}^d)}\|\Lambda^su\|_{L^{p_4}(\mathbb{T}^d)}),$$
with $p_i\in (1,\infty], i=1,...,4$ such that
$$\frac{1}{p}=\frac{1}{p_1}+\frac{1}{p_2}=\frac{1}{p_3}+\frac{1}{p_4}.$$
\vskip.10in
\textbf{Estimates on the whole space}

We also collect some important properties for the weighted Besov spaces from [MW15] and [Tri06], which are parallel to those for Besov spaces on the torus.
The  Schauder estimate takes the following form:

\vskip.10in
\th{Lemma 2.5}([MW15, Propositions 3.11, 3.12]) (i) Let $u\in \hat{\mathcal{B}}^{\alpha,\sigma}_{p,q}$ for some $\alpha\in \mathbb{R}$, $p,q\in[1,\infty]$. Then we have for every $\delta\geq0$
$$\|e^{tA}u\|_{\hat{\mathcal{B}}^{\alpha+\delta,\sigma}_{p,q}}\lesssim t^{-\delta/2}\|u\|_{\hat{\mathcal{B}}^{\alpha,\sigma}_{p,q}}.$$
(ii) Let $\alpha\leq \beta\in\mathbb{R}$ be such that $\beta-\alpha\leq 2$, $\sigma>0$ and $p, q\in[1,\infty]$. Then for $u\in\hat{\mathcal{B}}^{\beta,\sigma}_{p,q}$
$$\|(1-e^{tA})u\|_{\hat{\mathcal{B}}^{\alpha,\sigma}_{p,q}}\lesssim t^{\frac{\beta-\alpha}{2}}\|u\|_{\hat{\mathcal{B}}^{\beta,\sigma}_{p,q}}.$$

\vskip.10in
The multiplicative structure and Besov embedding theorems can be written as follows:

\vskip.10in
\th{Lemma 2.6}([MW15, Corollary 3.19, Corollary 3.21]) (1) For $\alpha>0, p_1, p_2, p, q\in[1,\infty],  \frac{1}{p_1}+\frac{1}{p_2}=\frac{1}{p}$, the bilinear map $(u; v)\mapsto uv$
extends to a continuous map from $\hat{\mathcal{B}}^{\alpha,\sigma}_{p_1,q}\times \hat{\mathcal{B}}^{\alpha,\sigma}_{p_2,q}$ to $\hat{\mathcal{B}}^{\alpha,\sigma}_{p,q}$.

(2) For $\alpha<0, \alpha+\beta>0,  p_1, p_2, p, q\in[1,\infty],   \frac{1}{p_1}+\frac{1}{p_2}=\frac{1}{p}$, the bilinear map $(u; v)\mapsto uv$
extends to a continuous map from $\hat{\mathcal{B}}^{\alpha,\sigma}_{p_1,q}\times \hat{\mathcal{B}}^{\beta,\sigma}_{p_2,q}$ to $\hat{\mathcal{B}}^{\alpha,\sigma}_{p,q}$.

(3) (Besov embedding [Tri06, Chapter 6]) Let $\alpha_1\leq \alpha_2$, $1\leq p_1\leq p_2\leq \infty$,
and $1\leq q_1\leq q_2\leq \infty$. Then
$$\hat{\mathcal{B}}^{\alpha_2,\sigma}_{p_1,q_1}\subset \hat{\mathcal{B}}^{\alpha_1,\sigma}_{p_1,q_1};\quad \hat{\mathcal{B}}^{\alpha_1,\sigma}_{p_1,q_1}\subset \hat{\mathcal{B}}^{\alpha_1,\sigma}_{p_1,q_2}.$$
If $\sigma>d$, then
$$\hat{\mathcal{B}}^{\alpha_1,\sigma}_{p_2,q_1}\subset \hat{\mathcal{B}}^{\alpha_1,\sigma}_{p_1,q_1}.$$
Here  $\subset$ means that the embedding is continuous and dense.
\vskip.10in

\section{Finite volume case}
In this section we consider (1.1) on the torus $\mathbb{T}^2$.

\subsection{Wick power}
In the following we define the Wick powers. First we define Wick powers on $L^2(\mathcal{S}'(\mathbb{T}^2),\mu)$ with $\mu=\frac{1}{2}N(0,(-\Delta+1)^{-1}):=N(0,C)$.
 \vskip.10in
\textbf{Wick powers on $L^2(\mathcal{S}'(\mathbb{T}^2),\mu)$}

In fact $\mu$ is a measure supported on $\mathcal{S}'(\mathbb{T}^2)$. We have the well-known (Wiener-It\^{o}) chaos decomposition $$L^2(\mathcal{S}'(\mathbb{T}^2),\mu)=\bigoplus_{n\geq0}\mathcal{H}_n.$$
Now we define the Wick powers by using approximations: for $\phi\in \mathcal{S}'(\mathbb{T}^2)$ define $$\phi_\varepsilon:=\rho_\varepsilon*\phi,$$ with $\rho_\varepsilon$ an approximate delta function on $\mathbb{R}^2$ given by
$$\rho_\varepsilon(x)=\varepsilon^{-2}\rho(\frac{x}{\varepsilon})\in \mathcal{D}, \int \rho=1.$$
Here the convolution means that we view $\phi$ as a periodic distribution in $\mathcal{S}'(\mathbb{R}^2)$ and do convolution on $\mathbb{R}^2$.
 For every $n\in\mathbb{N}$ we set $$:\phi_\varepsilon^n:_C:=c_\varepsilon^{n/2}P_n(c_\varepsilon^{-1/2}\phi_\varepsilon),$$
where $P_n,n=0,1,...,$ are the Hermite polynomials defined by the formula
$$P_n(x)=\sum_{j=0}^{[n/2]}(-1)^j\frac{n!}{(n-2j)!j!2^j}x^{n-2j},$$
and
$c_\varepsilon=\int\phi^2_\varepsilon\mu(d\phi)=\int\int \bar{G}(x-y)\rho_\varepsilon(y)dy\rho_\varepsilon(x)dx
=\|\bar{K}_\varepsilon\|_{L^2(\mathbb{R}\times \mathbb{T}^2 )}^2$. Then $$:\phi_\varepsilon^n:_C\in \mathcal{H}_n.$$
Here and in the following $\bar{G}$ is the Green function associated with $-A$ on $\mathbb{T}^2$ and  $\bar{K}(t,x)$ is the heat kernel associated with $A$ on $\mathbb{T}^2$ and $\bar{K}_\varepsilon=\bar{K}*\rho_\varepsilon$, where $*$ means  convolution in space and we view $\bar{K}$ as a periodic function on $\mathbb{R}^2$.

For Hermite polynomial $P_n$ we have that for $s,t\in\mathbb{R}$
\begin{equation}P_n(s+t)=\sum_{m=0}^nC_n^mP_m(s)t^{n-m},\end{equation}
where $C_n^m=\frac{n!}{m!(n-m)!}$.

 \vskip.10in
A direct calculation yields the following:
 \vskip.10in

\th{Lemma 3.1} Let $\alpha<0$, $n\in\mathbb{N}$ and $p>1$. $:\phi_\varepsilon^n:_C$ converges to some element in $L^p(\mathcal{S}'(\mathbb{T}^2),\mu;\mathcal{C}^{\alpha})$. This limit is called the $n$-th Wick power of $\phi$ with respect to the covariance $C$ and denoted by $:\phi^n:_C$.

\proof In fact,
  for every $p>1,\varepsilon_1,\varepsilon_2>0$, $m\in\mathbb{N}$ by (2.2) and (2.4) we have that
$$\aligned &\int\|:\phi_{\varepsilon_1}^m:_C-:\phi_{\varepsilon_2}^m:_C\|^{2p}_{\alpha}\mu(d\phi)\\\lesssim& \sum_{\psi\in{\Psi}_\star}\sum_{n\geq0}\sum_{x\in\Lambda_n}2^{2\alpha pn+2np}\int|\langle :\phi_{\varepsilon_1}^m:_C-:\phi_{\varepsilon_2}^m:_C,\psi_{x}^{n}\rangle|^{2p}\mu(d\phi)w(x)^{2p}
\\\lesssim&\sum_{\psi\in{\Psi}_\star}
\sum_{n\geq0}\sum_{x\in\Lambda_n}2^{2\alpha pn+2np}w(x)^{2p}
(\int|\langle :\phi_{\varepsilon_1}^m:_C-:\phi_{\varepsilon_2}^m:_C,\psi_{x}^{n}\rangle|^{2}\mu(d\phi))^p,\endaligned$$
where $\sigma>0$ in $w(x)$ and  in the last inequality we used the hypercontractivity of the Gaussian measure. Moreover, we obtain that
$$\aligned&\int|\langle:\phi_{\varepsilon_1}^m:_C-:\phi_{\varepsilon_2}^m:_C,\psi_{x}^{n}\rangle|^{2}\mu(d\phi)
\\\lesssim&
 \int\int |\psi_{x}^{n}(y)\psi_{x}^{n}(\bar{y})|\bigg{|}(\int\phi_{\varepsilon_1}(y)\phi_{\varepsilon_1}(\bar{y})\mu(d\phi))^m
 -2(\int\phi_{\varepsilon_1}(y)\phi_{\varepsilon_2}(\bar{y})\mu(d\phi))^m\\&+(\int\phi_{\varepsilon_2}(y)\phi_{\varepsilon_2}(\bar{y})\mu(d\phi))^m\bigg{|}
 dyd\bar{y}
 \\\lesssim&
 \int\int |\psi_{x}^{n}(y)\psi_{x}^{n}(\bar{y})|\bigg{|}(\int\int\rho_{\varepsilon_1}(y-x_1)\rho_{\varepsilon_1}(\bar{y}-x_2)\bar{G}(x_1-x_2)dx_1dx_2)^m
 \\&-2(\int\int\rho_{\varepsilon_1}(y-x_1)\rho_{\varepsilon_2}(\bar{y}-x_2)\bar{G}(x_1-x_2)dx_1dx_2)^m
 \\&+(\int\int\rho_{\varepsilon_2}(y-x_1)\rho_{\varepsilon_2}(\bar{y}-x_2)\bar{G}(x_1-x_2)dx_1dx_2)^m\bigg{|}dyd\bar{y}
   \\\lesssim&
(\varepsilon_1^\kappa+\varepsilon_2^\kappa) \int\int |\psi_{x}^{n}(y)\psi_{x}^{n}(\bar{y})| |y-\bar{y}|^{-\delta}dyd\bar{y}\lesssim (\varepsilon_1^\kappa+\varepsilon_2^\kappa)2^{-2n+n\delta},\endaligned$$
 where $\delta>\kappa>0,2\alpha+\delta<0$. Here in the third inequality we have used Lemma 10.17 in [Hai14]. Thus the results follow from a direct calculation.$\hfill\Box$
 \vskip.10in

\textbf{Wick powers on a fixed probability space}

Now  we follow the idea from [DD03] and [MW15] to define the Wick powers of the solutions to the stochastic heat equation in the paths space.
We fix a probability space $(\Omega,\mathcal{F},P)$ and $W$ is a cylindrical Wiener process on $L^2(\mathbb{T}^2)$. We also have the well-known (Wiener-It\^{o}) chaos decomposition $$L^2(\Omega,\mathcal{F},P)=\bigoplus_{n\geq0}\mathcal{H}'_n.$$In the following we set $Z(t)=\int_0^t e^{(t-s)A}dW(s)$, and we can also define Wick powers of $Z(t)$ with respect to different covariances by approximations: Let $Z_\varepsilon(t,x)=\int_0^t \langle\bar{K}_\varepsilon(t-s,x-\cdot),dW(s)\rangle$. Here $\langle\cdot,\cdot\rangle$ means  inner product in $L^2(\mathbb{T}^2)$.  For every $n\in\mathbb{N}$ we set $$:Z_\varepsilon^n(t):_{C_t}:=(c_{\varepsilon,t})^{\frac{n}{2}}P_n((c_{\varepsilon,t})^{-\frac{1}{2}}Z_\varepsilon(t))\in \mathcal{H}'_n,$$
where $P_n,n=0,1,...,$ are the Hermite polynomials
and $c_{\varepsilon,t}=\|1_{[0,t]}\bar{K}_\varepsilon\|_{L^2(\mathbb{R}\times \mathbb{T}^2)}^2$.
\vskip.10in

In the following we prove that $:Z_\varepsilon^n(t):_{C_t}$ is a Cauchy sequence  in $C([0,T];\mathcal{C}^\alpha)$ and  define the Wick powers of $Z$ as the limit. To prove this we have to use the following result from [ZZ15, Lemma 4.1], the proof of which is a modification of the proof of [Hai14, Lemma 10.18].
\vskip.10in

\th{Lemma 3.2} If $|K(t,y)|\lesssim (t+|y|^2)^{\zeta/2}$ with $\zeta\in (-3,0)$, then
$$|(K*\rho_{\varepsilon})(t,y)|\leq Ct^{-\frac{\delta}{2}}(t+|y|^2)^{\frac{\zeta+\delta}{2}},$$
for $0<\delta<1, \zeta+\delta>-2$,
and
$$|(K*\rho_{\varepsilon})(t,y)-K(t,y)|\lesssim t^{-\frac{\delta}{2}}\varepsilon^{\zeta-\bar{\zeta}}|y|^{\bar{\zeta}+\delta},$$
for $\bar{\zeta}+\delta>-2$.

 \vskip.10in
Using Lemma 3.2 we obtain the Wick powers of $Z(t)$.
 \vskip.10in

 \th{Lemma 3.3} For  $\alpha<0$, $n\in\mathbb{N}$, $p>1$,  $:Z_\varepsilon^n(t):_{C_t}$ converges  in $L^p(\Omega,C([0,T];\mathcal{C}^{\alpha}))$. The limit is called Wick power of $Z(t)$ with respect to the covariance $C_t$ and denoted by $:Z^n(t):_{C_t}$.

\proof We first prove that $Z_\varepsilon\in C([0,T];\mathcal{C}^\alpha)$ $P$-almost-surely. By the factorization method in [D04] we have that for $\kappa\in(0,1)$
$$Z_\varepsilon(t)=\frac{\sin(\pi\kappa)}{\pi}\int_0^t(t-s)^{\kappa-1}\langle\bar{K}_\varepsilon(t-s,x-\cdot),U(s)\rangle ds,$$
where $$U(s,y)=\int_0^s(s-r)^{-\kappa} \langle\bar{K}(s-r,y-\cdot),dW(r)\rangle.$$
A similar argument as in the proof of Lemma 2.7 in [D04] implies that it suffices to prove that for $p>1/(2\kappa),$ \begin{equation}\textbf{E}\|U\|_{L^{2p}(0,T;\mathcal{C}^\alpha)}<\infty.\end{equation}
In fact, by (2.2) and (2.4) we have that
$$\aligned \textbf{E}\|U(s)\|_\alpha^{2p}\lesssim &\sum_{\psi\in{\Psi}_\star}\sum_{n\geq0}\sum_{x\in\Lambda_n}\textbf{E}2^{2\alpha pn+2np}
|\langle U(s)
,\psi_{x}^{n}\rangle|^{2p}w(x)^{2p}\\\lesssim&\sum_{\psi\in{\Psi}_\star}\sum_{n\geq0}\sum_{x\in\Lambda_n}2^{2\alpha pn+2np}
(\textbf{E}|\langle U(s)
,\psi_{x}^{n}\rangle|^{2})^pw(x)^{2p}.\endaligned$$
Here $\sigma>0$ in $w(x)$ and we used Gaussian hypercontractivity in the second inequality.
Moreover we obtain that
$$\aligned \textbf{E}|\langle U(s),\psi_{x}^{n}\rangle|^{2}\leq&\int\int|\psi_{x}^{n}(y)\psi_{x}^{n}(\bar{y})|\int_0^s(s-r)^{-2\kappa}\bar{K}*\bar{K}(s-r,y-\bar{y})drdyd\bar{y}
\\\lesssim&\int\int|\psi_{x}^{n}(y)\psi_{x}^{n}(\bar{y})|\int_0^s(s-r)^{\kappa-1}|y-\bar{y}|^{-8\kappa}drdyd\bar{y}
\\\lesssim&2^{-2n+8n\kappa}s^\kappa,\endaligned$$
where we used [Hai14, Lemma 10.17] in the second inequality. Thus (3.2) follows by choosing $\kappa$ small enough and a direct calculation.
 Now we prove that $:Z_{\varepsilon}^m:$ is a Cauchy sequence.
  For every $p>1$, by (2.2) and (2.4) we have  for $t_1,t_2\geq0$ that
$$\aligned &\textbf{E}\|(:Z_{\varepsilon_1}^m:_{C_{t_1}}-:Z^m_{\varepsilon_2}:_{C_{t_1}})(t_1,\cdot)-(:Z_{\varepsilon_1}^m:_{C_{t_2}}
-:Z^m_{\varepsilon_2}:_{C_{t_2}})(t_2,\cdot)\|^{2p}_{\alpha}\\\leq& \sum_{\psi\in{\Psi}_\star}\sum_{n\geq0}\sum_{x\in\Lambda_n}\textbf{E}2^{2\alpha pn+2np}
|\langle (:Z_{\varepsilon_1}^m:_{C_{t_1}}-:Z^m_{\varepsilon_2}:_{C_{t_1}})(t_1,\cdot)-(:Z_{\varepsilon_1}^m:_{C_{t_2}}
-:Z^m_{\varepsilon_2}:_{C_{t_2}})(t_2,\cdot)
,\psi_{x}^{n}\rangle|^{2p}w(x)^{2p}
\\\lesssim&\sum_{\psi\in{\Psi}_\star}
\sum_{n\geq0}\sum_{x\in\Lambda_n}2^{2\alpha pn+2np}
(\textbf{E}|\langle (:Z_{\varepsilon_1}^m:_{C_{t_1}}-:Z^m_{\varepsilon_2}:_{C_{t_1}})(t_1,\cdot)-(:Z_{\varepsilon_1}^m:_{C_{t_2}}
-:Z^m_{\varepsilon_2}:_{C_{t_2}})(t_2,\cdot),\psi_{x}^{n}\rangle|^{2})^pw(x)^{2p},\endaligned$$
where we used Gaussian hypercontractivity in the second inequality.
For  convenience we use $\xi$ to denote space-time white noise given by $\int \phi(s,y)\xi(ds,dy)=\int_{\mathbb{R}^+}\langle\phi, dW(s)\rangle$ for $\phi\in L^2(\mathbb{R}^+\times \mathbb{T}^2)$.
Then we obtain that for $k=1,2$ and $j=1,2$ $$:Z_{\varepsilon_k}^m(t_j):_{C_{t_j}}=\int\Pi_{i=1}^m\bar{K}_{\varepsilon_k}({t_j-s_i},y-y_i)1_{s_i\in[0,t_1]}\xi(d\eta_1)...\xi(d\eta_m),$$
 where $\eta_a=(s_a,y_a)$,
and $\int f(\eta_{1...n})\xi(d\eta_1)...\xi(d\eta_m)$ denotes
  a generic element of the $n$-th chaos of $\xi$ for $\eta_{1...n}=\eta_1...\eta_n$.
Moreover, for $t_1\leq t_2$ to estimate
 $$\textbf{E}|\langle(:Z_{\varepsilon_1}^m:_{C_{t_1}}-:Z^m_{\varepsilon_2}:_{C_{t_1}})(t_1,\cdot)-(:Z_{\varepsilon_1}^m:_{C_{t_2}}
-:Z^m_{\varepsilon_2}:_{C_{t_2}})(t_2,\cdot),\psi_{x}^{n}\rangle|^{2},$$ it suffices to calculate
$$\aligned&\int|\langle\Pi_{i=1}^m\bar{K}_{\varepsilon_1}({t_1-s_i},\cdot-y_i)1_{s_i\in[0,t_1]}
-\Pi_{i=1}^m\bar{K}_{\varepsilon_2}({t_1-s_i},\cdot-y_i)1_{s_i\in[0,t_1]}
\\&-[\Pi_{i=1}^m\bar{K}_{\varepsilon_1}({t_2-s_i},\cdot-y_i)1_{s_i\in[0,t_2]}
-\Pi_{i=1}^m\bar{K}_{\varepsilon_2}({t_2-s_i},\cdot-y_i)1_{s_i\in[0,t_2]}],\psi_{x}^{n}\rangle|^2d\eta_{1...m},\endaligned$$
which is bounded by
$$\aligned&2\int|\langle(\Pi_{i=1}^m\bar{K}_{\varepsilon_1}({t_1-s_i},\cdot-y_i)-\Pi_{i=1}^m\bar{K}_{\varepsilon_1}({t_2-s_i},\cdot-y_i))1_{s_i\in[0,t_1]}
\\&-(\Pi_{i=1}^m\bar{K}_{\varepsilon_2}({t_1-s_i},\cdot-y_i)-\Pi_{i=1}^m\bar{K}_{\varepsilon_2}({t_2-s_i},\cdot-y_i))1_{s_i\in[0,t_1]},\psi_{x}^{n}\rangle|^2d\eta_{1...m}
\\&+2\int|\langle[\Pi_{i=1}^m\bar{K}_{\varepsilon_1}({t_2-s_i},\cdot-y_i)1_{s_i\in[t_1,t_2]}
-\Pi_{i=1}^m\bar{K}_{\varepsilon_2}({t_2-s_i},\cdot-y_i)1_{s_i\in[t_1,t_2]}],\psi_{x}^{n}\rangle|^2d\eta_{1...m}.\endaligned$$
Since by Lemma 3.2 and [Hai14, Lemma 10.18]
we have that for every $\delta>0$
$$\aligned&|\bar{K}_{\varepsilon_1}({t_1-s_i},y-y_i)-\bar{K}_{\varepsilon_1}({t_2-s_i},y-y_i)|\\\lesssim& |t_1-t_2|^\delta (|t_2-s_i|^{-\frac{1}{2}+\frac{\delta}{2}}+|t_1-s_i|^{-\frac{1}{2}+\frac{\delta}{2}})|y-y_i|^{-1-3\delta},\endaligned$$
and $$\aligned&|\bar{K}_{\varepsilon_1}({t_1-s_i},y-y_i)-\bar{K}_{\varepsilon_2}({t_1-s_i},y-y_i)|\\\lesssim& (\varepsilon_1^{2\delta}+\varepsilon_2^{2\delta}) |t_1-s_i|^{-\frac{1}{2}+\frac{\delta}{2}}|y-y_i|^{-1-3\delta},\endaligned$$
which combining with the interpolation and [Hai14, Lemma 10.14] implies that for every $\delta>0$
$$\aligned& \textbf{E}|\langle(:Z_{\varepsilon_1}^m:_{C_{t_1}}-:Z^m_{\varepsilon_2}:_{C_{t_1}})(t_1,\cdot)-(:Z_{\varepsilon_1}^m:_{C_{t_2}}
-:Z^m_{\varepsilon_2}:_{C_{t_2}})(t_2,\cdot),\psi_{x}^{n}\rangle|^{2}\\\lesssim&(\varepsilon_1^{2\delta}+\varepsilon_2^{2\delta})|t_2-t_1|^{\delta}\int\int |\psi^{n}(y)\psi^{n}(\bar{y})||y-\bar{y}|^{-6\delta} dyd\bar{y}\\\lesssim& (\varepsilon_1^{2\delta}+\varepsilon_2^{2\delta})|t_2-t_1|^{\delta}2^{-2n+6n\delta}.\endaligned$$
Then the above estimates yield that $$\aligned &\textbf{E}\|(:Z_{\varepsilon_1}^m:_{C_{t_1}}-:Z^m_{\varepsilon_2}:_{C_{t_1}})(t_1,\cdot)-(:Z_{\varepsilon_1}^m:_{C_{t_2}}
-:Z^m_{\varepsilon_2}:_{C_{t_2}})(t_2,\cdot)\|^{2p}_{\alpha}\\\lesssim&\sum_{\psi\in\bar{\Psi}_\star}
\sum_{n\geq0}2^{2\alpha np+2np+2n}
(\varepsilon_1^{2\delta}+\varepsilon_2^{2\delta})^{ p}|t_2-t_1|^{\delta p}2^{-2np+6np\delta}.\endaligned$$
Thus the results follow from Kolmogorov's continuity test (in time) if we choose $\delta>0$ small enough and $p$ sufficiently large. $\hfill\Box$
 \vskip.10in
\th{Remark} We can also use the approximations from [DD03] to define the
Wick powers. We can prove that these two approximations converge to the same limit. This follows from the fact that  under $\mu$, $\phi=^d\int_{\mathbb{T}^2}\int_{-\infty}^t\bar{K}(t-s,\cdot-y)\xi(ds,dy)$ and for every $\varphi\in \mathcal{S}(\mathbb{T}^2)$, $$\langle:\phi^n:,\varphi\rangle=^d\int_{[(-\infty,t]\times \mathbb{T}^2]^n}\langle\varphi,\prod_{j=1}^n\bar{K}(t-r_j,\cdot-y_j)\rangle\xi(dr_1,dy_1)...\xi(dr_n,dy_n) $$
Here  $\xi$ is  space-time white noise.

 \vskip.10in

By this lemma we can also define the Wick powers with respect to another covariance $:Z^n(t):_C$.
 \vskip.10in

\th{Lemma 3.4} For  $\alpha<0$, $p>1$,  $n\in\mathbb{N}$ and $t>0$, $:Z_\varepsilon^n(t):_{C}:=c_\varepsilon^{\frac{n}{2}}P_n(c_\varepsilon^{-\frac{1}{2}}Z_\varepsilon (t))$ converges  in $L^p(\Omega,C((0,T];\mathcal{C}^{\alpha}))$. Here the norm for $C((0,T];\mathcal{C}^{\alpha})$ is $\sup_{t\in[0,T]}t^{\rho/2}\|\cdot\|_\alpha$ for  $\rho>0$. The limit is called Wick powers of $Z(t)$ with respect to the covariance $C$ and denoted by $:Z^n(t):_C$. Moreover, for $t>0$
$$:Z^n(t):_C=\sum_{l=0}^{[n/2]}c_t^l\frac{n!}{(n-2l)!l!2^l}:Z^{n-2l}(t):_{C_t},$$
where $c_t:=\lim_{\varepsilon\rightarrow0}(c_{\varepsilon,t}-c_\varepsilon)$ locally uniformly for $t\in(0,T]$.

\proof By Lemma 3.3 it follows that for every $n\in \mathbb{N}$, $p>1$,
$$:Z_\varepsilon^n(t):_{C_t}\rightarrow :Z^n(t):_{C_t}\quad \textrm{ in } L^p(\Omega,C([0,T];\mathcal{C}^{\alpha})).$$
By the definition of $c_{\varepsilon,t}$ and $c_\varepsilon$ we also have that for every $\rho>0$, $t>0$ and $\varepsilon>0$
$$|c_{\varepsilon,t}-c_\varepsilon|\lesssim t^{-\rho/2},$$
and
$$c_t:=\lim_{\varepsilon\rightarrow0}(c_{\varepsilon,t}-c_\varepsilon)=-\int_t^\infty\int_{\mathbb{T}^2} \bar{K}(r,x)^2dxdr.$$

Moreover, the definition of $P_n$ yields that
$$\aligned :Z_\varepsilon^n(t):_C
=\sum_{l=0}^{[n/2]}(c_{\varepsilon,t}-c_\varepsilon)^l\frac{n!}{(n-2l)!l!2^l}:Z_\varepsilon^{n-2l}(t):_{C_t},
\endaligned$$
which implies the result by letting $\varepsilon\rightarrow0$.$\hfill\Box$

\vskip.10in

Now following the technique in [MW15] we combine the initial value part  with the Wick powers by using (3.1).
We set $V(t)=e^{tA}z,V_\varepsilon=\rho_\varepsilon*V$, $z\in{\mathcal{C}}^{\alpha}$ for $\alpha<0$ and
$$\bar{Z}(t)=Z(t)+V(t),\quad \bar{Z}_\varepsilon(t)=Z_\varepsilon(t)+V_\varepsilon(t),$$
$$:\bar{Z}^n(t):_C=\sum_{k=0}^nC_n^kV(t)^{n-k}:Z^k(t):_C,\quad :\bar{Z}_\varepsilon^n(t):_C=\sum_{k=0}^nC_n^kV_\varepsilon(t)^{n-k}:Z_\varepsilon^k(t):_C.$$
By Lemma 2.2 we know that $V\in C([0,T],\mathcal{C}^\alpha)$  and  $V\in C((0,T],\mathcal{C}^\beta)$ for $\beta>\alpha$ with the norm $\sup_{t\in[0,T]}t^{\frac{\beta-\alpha}{2}}\|\cdot\|_{\beta}$. Moreover, $$\sup_{t\in[0,T]}t^{\frac{\beta-\alpha}{2}}\|V(t)\|_{\beta}\lesssim \|z\|_{\alpha},$$ for $\beta>\alpha$. Then by Lemmas 2.3 and 3.4 we have the following results:
\vskip.10in

\th{Lemma 3.5} Let $\alpha<0, z\in\mathcal{C}^\alpha$, $p>1$. Define $\bar{Z}$ and $:\bar{Z}^n:_C$ as above. Then for   $n\in\mathbb{N}$  $:\bar{Z}^n_\varepsilon:_C$ converges to $:\bar{Z}^n:_C$ in  $L^p(\Omega,C((0,T]; {\mathcal{C}}^{\alpha}))$.  Here the norm for $C((0,T];\mathcal{C}^\alpha)$ is $$\sup_{t\in[0,T]}t^\frac{(\beta-\alpha)n+\rho}{2}\|\cdot\|_{{\alpha}}$$ for $\beta+\alpha>0, \rho>0$.

\vskip.10in
\textbf{Relations between two different Wick powers}

  First we introduce the following measure. Set $:q(\phi):=\sum_{n=0}^{2N}a_n:\phi^n:_C$,   $:p(\phi):=\sum_{n=1}^{2N}na_n:\phi^{n-1}:_C$ and we assume that $a_n\in \mathbb{R}$ and $a_{2N}>0$. Let $$\nu=c\exp{(-\int_{\mathbb{T}^2}:q(\phi):dx)}\mu,$$
where $c$ is a normalization constant. Then by  [GlJ86, Sect. 8.6] for every $p\in [1,\infty)$, $\varphi(\phi):=\exp{(-\int_{\mathbb{T}^2}:q(\phi):dx)}\in L^p(\mathcal{S}'(\mathbb{T}^2),\mu)$. The following result states the relations between two different Wick powers.
\vskip.10in

\th{Lemma 3.6} Let $\phi$ be a measurable map from $(\Omega, \mathcal{F},P)$ to $C([0,T],B^{-\gamma}_{2,2})$ with $\gamma>2$, $P\circ \phi(t)^{-1}=\nu$ for every $t\in[0,T]$ and let $\bar{Z}(t)$ be defined as above. Assume in addition that  $y=\phi-\bar{Z}\in  C([0,T];\mathcal{C}^{\beta})$ $P$-a.s. for some $\beta>-\alpha>0$. Then for every $t>0$, $n\in\mathbb{N}$
$$:\phi^n(t):_C=\sum_{k=0}^nC_n^ky^{n-k}(t):\bar{Z}^{k}(t):_{C}\quad P-a.s..$$

\proof By Lemma 3.5 it follows that for every $k\in \mathbb{N}$, $p>1$
$$:\bar{Z}_\varepsilon^k:_{C}\rightarrow :\bar{Z}^k:_{C}\quad \textrm{ in } L^p(\Omega, C((0,T];\mathcal{C}^{\alpha})), \textrm{ as }\varepsilon\rightarrow0.$$
Since $y_\varepsilon=\phi_\varepsilon-\bar{Z}_{\varepsilon}=\rho_\varepsilon*y$ and $y\in  C([0,T];\mathcal{C}^{\beta})$ $P$-a.s., it is obvious that  $y_\varepsilon\rightarrow y$ in $C([0,T];\mathcal{C}^{\beta-\kappa})$ $P$-a.s. for every $\kappa>0$ with $\beta-\kappa+\alpha>0$, which combined with Lemma 2.3 implies that for $k\in\mathbb{N}$, $k\leq n$,
$$y_\varepsilon^{n-k}:\bar{Z}_\varepsilon^{k}:_{C}\rightarrow^P y^{n-k}:\bar{Z}^{k}:_{C} \quad \textrm{in } C((0,T];\mathcal{C}^{\alpha}),  \textrm{ as }\varepsilon\rightarrow0.$$
Since $\exp{(-\int_{\mathbb{T}^2}:q(\phi):dx)}\in L^p(\mathcal{S}'(\mathbb{T}^2),\mu)$ for every $p\geq 1$, by H\"{o}lder's inequality and Lemma 3.1 we get that for $t>0$ and $p>1$
$$:\phi_\varepsilon^n(t):_C\rightarrow :\phi^n(t):_C \quad \textrm{ in } L^p(\Omega,\mathcal{C}^{\alpha}),  \textrm{ as }\varepsilon\rightarrow0.$$
Moreover, by (3.1) we have
$$\aligned &:\phi_\varepsilon^n:_C=:(y_\varepsilon+\bar{Z}_\varepsilon)^n:_C=c_\varepsilon^{n/2}P_n(c_\varepsilon^{-1/2}(y_\varepsilon+\bar{Z}_\varepsilon))
\\=&\sum_{k=0}^nC_n^k c_\varepsilon^{n/2}P_k(c_\varepsilon^{-1/2}\bar{Z}_\varepsilon)(c_\varepsilon^{-1/2}y_\varepsilon)^{n-k}\\=&\sum_{k=0}^nC_n^k:\bar{Z}_\varepsilon^k:_C y_\varepsilon^{n-k},
\endaligned$$
which implies the result by letting $\varepsilon\rightarrow0$.
$\hfill\Box$
\vskip.10in

In the following, we only use Wick powers $:\cdot:_C$ and we write  $:\cdot:$ for simplicity.

\subsection{Relations between the two solutions: starting with solutions given by Dirichlet forms}
As mentioned in the introduction,   weak solutions to (1.1) have been constructed in [AR91] by Dirichlet forms. In this subsection we prove that the solutions constructed in [AR91] also satisfy the shifted equation. First we recall some basic results related to Dirichlet forms from [AR91].
\vskip.10in

{\textbf{Solutions given by Dirichlet forms}}

Let $H=L^2(\mathbb{T}^2)$ and let $-\Delta + I$ be the generator of the
following quadratic form on $H : (u,v)\mapsto \int_{\mathbb{T}^2} \langle\nabla u, \nabla v\rangle_{\mathbb{R}^d} dx + \int_{\mathbb{T}^2} uv dx$
with $u, v\in \{g\in L^2(\mathbb{T}^2)|\nabla g\in L^2(\mathbb{T}^2)\}$ (where $\nabla$ is in the sense of
distributions). Let $\{e_k| k\in\mathbb{Z}^2\}\subset C^\infty({\mathbb{T}^2})$ be the (orthonormal) eigenbasis of
$-\Delta+I$ in $H$ and $\{\lambda_k|k\in\mathbb{Z}^2\}\subset (0,\infty)$ the corresponding eigenvalues. Define for $s\in\mathbb{R}$,
$$H^s:=\{u\in \mathcal{S}'(\mathbb{T}^2) |\sum_{k\in\mathbb{Z}^2}\lambda_k^s{ }_{\mathcal{S}'}\langle u,e_k\rangle^2_{\mathcal{S}}<\infty\},$$
equipped with the inner product
$$\langle u,v\rangle_{H^s}:=\sum_{k\in\mathbb{Z}^2}\lambda_k^s{ }_{\mathcal{S}'}\langle u,e_k\rangle_{\mathcal{S}}{ }_{\mathcal{S}'}\langle v,e_k\rangle_{\mathcal{S}}.$$
 If for $s\geq0$ ${ }_{H^s}\langle\cdot,\cdot\rangle_{H^{-s}}$ denotes  the dualization
between  $H^s$ and its dual space $H^{-s}$,  then it follows that
$$ { }_{H^s}\langle u, v\rangle_{H^{-s}}=\langle u, v\rangle_{H},   u\in H^s ,v\in {H}.$$
Let $E=H^{-1-\epsilon}, E^*=H^{1+\epsilon}$ for some $\epsilon>0$.  We denote their Borel $\sigma$-algebras by $\mathcal{B}(E), \mathcal{B}(E^*)$ respectively.
Define $$\mathcal{F}C_b^\infty=\{u:u(z)=f({ }_{E^*}\!\langle l_1,z\rangle_E,{ }_{E^*}\!\langle l_2,z\rangle_E,...,{ }_{E^*}\!\langle l_m,z\rangle_E),z\in E, l_1,l_2,...,l_m\in E^*, m\in \mathbb{N}, f\in C_b^\infty(\mathbb{R}^m)\}.$$
 Define for $u\in \mathcal{F}C_b^\infty$ and $l\in H$, $$\frac{\partial u}{\partial l}(z):=\frac{d}{ds}u(z+sl)|_{s=0},z\in E,$$
  that is, by the chain rule,
  $$\frac{\partial u}{\partial l}(z)=\sum_{j=1}^m\partial_jf({ }_{E^*}\!\langle l_1,z\rangle_E,{ }_{E^*}\!\langle l_2,z\rangle_E,...,{ }_{E^*}\!\langle l_m,z\rangle_E)\langle l_j,l\rangle.$$
Let $Du$ denote the $H$-derivative of $u\in \mathcal{F}C_b^\infty$, i.e. the map from $E$ to $H$ such that $$\langle Du(z),l\rangle=\frac{\partial u}{\partial l}(z)\textrm{ for all } l\in H, z\in E.$$
By [AR91] we easily deduce that the form
$$\mathcal{E}(u,v):=\frac{1}{2}\int_E \langle Du, Dv\rangle_{H}d\nu; u,v\in\mathcal{F}C_b^\infty$$
is closable and its closure $(\mathcal{E},D(\mathcal{E}))$ is a quasi-regular Dirichlet form on $L^2(E;\nu)$ in the sense of [MR92]. By [AR91, Theorem 3.6] we know that there exists a (Markov) diffusion process $M=(\Omega,\mathcal{F},\mathcal{M}_t,(X(t))_{t\geq0},(P^z)_{z\in E})$ on $E$ \emph{properly associated with} $(\mathcal{E},D(\mathcal{E}))$, i.e. for $u\in L^2(E;\nu)\cap\mathcal{B}_b(E)$, the transition semigroup $P_tu(z):=E^z[u(X(t))]$ is $\mathcal{E}$-quasi-continuous for all $t >0$ and is a $\nu$-version of $T_tu$, where $T_t$ is the semigroup associated with $(\mathcal{E},D(\mathcal{E}))$. Here for the notion of $\mathcal{E}$-quasi-continuity we refer to [MR92, ChapterIII, Definition 3.2].
 \vskip.10in

By [GlJ86, (9.1.32)] we have the following:
 \vskip.10in

\th{Theorem3.7} For each $l$ smooth, we have that the partial  log derivative $\beta_l$ of $\nu$ is given by $$\beta_l(z)= -\sum_{n=1}^{2N}na_n:z^{n-1}:( l)+_{H^{s}}\!\langle \Delta l-l,z\rangle_{H^{-s}},$$
where $:z^{n}:( l)$ denotes the dualization between $:z^n:$ and $l$.

 \vskip.10in
Now we want to extend the definition of $\beta_l$ to the whole space $E$.  By [R86, Theorem 3.1] $l\rightarrow :z^n:(l)$ can be extended to a continuous map from $H$ to $ L^2(E,\mu)$. So, by [AR91, Proposition 6.9], there exists a $\mathcal{B}(H^{-1-\epsilon})/ \mathcal{B}(H^{-1-\epsilon})$ measurable map $:z^n::H^{-1-\epsilon}\rightarrow H^{-1-\epsilon}$ such that $:z^n:(l)=_{H^{-1-\epsilon}}\langle :z^n:,l\rangle_{H^{1+\epsilon}}$ $\nu$-a.e.. By [AR91, Theorem 6.10] we have the following Fukushima decomposition for $X(t)$ under $P^z$.
 \vskip.10in

\th{Theorem 3.8} There exist a map $W:\Omega\rightarrow C([0,\infty);E)$ and a \emph{properly  $\mathcal{E}$-exceptional set} $S\subset E$, i.e. $\nu(S)=0$ and $P^z[X(t)\in E\setminus S, \forall t\geq0]=1$ for $z\in E\backslash S$, such that $\forall z\in E\backslash S$ under $P^z$, $W$ is an $\mathcal{M}_t$- cylindrical Wiener process and the sample paths of the associated  process $M=(\Omega,\mathcal{F},(X(t))_{t\geq0},(P^z)_{z\in E})$ on $E$ satisfy the following:  for $l\in H^{2+s}$, $s>0$
 \begin{equation}\aligned{ }_{E^*}\!\langle l,X(t)-X(0)\rangle_{E}=&\int_0^t\langle l,dW(r)\rangle+\int_0^t\bigg[{ }_{H^{-s-2}}\!\langle -\sum_{n=1}^{2N}na_n:X(r)^{n-1}:,l\rangle_{H^{2+s}}\\&+_{H^{s}}\!\langle \Delta l-l,X(r)\rangle_{H^{-s}}\bigg] dr\quad \forall t\geq 0  \textrm{ }P^z\rm{-a.s.}.\endaligned\end{equation}
 Moreover, $\nu$ is an invariant measure for $M$ in the sense that $\int P_tu d\nu=\int ud\nu$ for $u\in L^2(E;\nu)\cap\mathcal{B}_b(E)$.
 \vskip.10in

\textbf{Relations between the two solutions}

 In the following we discuss the relations between $M$ constructed above and the shifted equation. In fact we have $\mathcal{C}^\alpha\subset E$ for $\alpha\in(-1,0)$, $\mathcal{C}^\alpha\in \mathcal{B}(E)$ and $\nu(\mathcal{C}^\alpha)=1$.
For $W$ constructed in Theorem 3.8 define
 $\bar{Z}(t):= \int_{0}^te^{(t-s)A}d W(s)+e^{tA}X(0)$.
\vskip.10in

\th{Theorem 3.9} Let $\alpha\in(-\frac{1}{2N-1},0)$, $-\alpha<\beta<\alpha+2$.  There exists a properly $\mathcal{E}$-exceptional set $S_2\subset E$ in the sense of Theorem 3.8  such that for every  $z\in\mathcal{C}^{\alpha}\setminus S_2$ under $P^z$,  $Y:=X-\bar{Z}\in C([0,T];\mathcal{C}^{\beta})$ is a solution to the following equation:
\begin{equation}Y(t)=-\int_0^te^{(t-s)A}\sum_{k=1}^{2N}ka_k\sum_{l=0}^{k-1} C_{k-1}^lY(s)^l:\bar{Z}(s)^{k-1-l}:ds.\end{equation}
Moreover,
$$P^z[X(t)\in {\mathcal{C}}^{\alpha}\setminus S_2, \forall t\geq0]=1 \textrm{ for } z\in \mathcal{C}^{\alpha}\setminus S_2.$$
 \proof Recall that $:p(\phi):=\sum_{n=1}^{2N}na_n:\phi^{n-1}:$.
 Now for $z\in E\setminus S$ under $P^z$
 we have that
$$X(t)=-\int_0^te^{(t-\tau)A}:p(X(\tau)):d\tau+\bar{Z}(t).$$
Since $\nu$ is an invariant measure for $X$, by Lemmas 2.1 and 3.1 we conclude that  for every $T\geq0$,  $p>1, \epsilon>0,$ with $ \alpha+2\epsilon<0,$ and $p_0>1$ large enough
$$\aligned&\int E^z\int_0^T\|:p(X(\tau)):\|_{{\alpha}}^pd\tau\nu(dz)\lesssim\int E^z\int_0^T\|:p(X(\tau)):\|_{B^{\alpha+\epsilon}_{p_0,p_0}}^pd\tau\nu(dz)
\\=& T\int\|:p(\phi):\|_{B^{\alpha+\epsilon}_{p_0,p_0}}^p\nu(d\phi)\lesssim T\int\|:p(\phi):\|_{\mathcal{C}^{\alpha+2\epsilon}}^p\nu(d\phi)<\infty,\endaligned$$
which implies that  there exists a properly $\mathcal{E}$-exceptional set $S_1\supset S$ such that for  $z\in E\setminus S_1$  $P^z$-a.s.
$$ :p(X(\cdot)):\in L^p(0,T;\mathcal{C}^\alpha),\quad  E^z\int_0^T\|:p(X(\tau)):\|_{{\alpha}}^pd\tau<\infty, \quad \forall p>1.$$
Here we used Lemma 2.1 to deduce the first result.  The second, however, does not imply the first  directly because of (2.1).
 Lemma 2.2 implies that for $\beta<\alpha+2$
$$\int_0^te^{(t-\tau)A}:p(X(\tau)):d\tau\in C([0,T];\mathcal{C}^{\beta})\quad P^z-a.s..$$
Now we conclude that for $z\in E\setminus S_1$
$$X-\bar{Z}\in  C([0,T];\mathcal{C}^{\beta})\quad P^z-a.s..$$
Since $P^\nu\circ X(t)^{-1}=\nu$, by Lemma 3.6 we conclude that under $P^\nu$, $Y:=X-\bar{Z}$ satisfies (3.4) and for $\nu$-a.e. $z\in E$ under $P^z$, $Y:=X-\bar{Z}$ satisfies (3.4).
In the following we prove that the results hold under $P^z$ for $z$ outside a properly $\mathcal{E}$-exceptional set.
Define $Z(t)=\int_0^te^{(t-s)A}dW(s)$ and we also obtain that
$$\aligned\bar{Y}(s,t_0):=&X(s+t_0)-Z(s+t_0)-e^{sA}(X(t_0)-Z(t_0))\\=&\int_{t_0}^{t_0+s}e^{(t_0+s-\tau)A}:p(X(\tau)):d\tau\in  C([0,T]^2;\mathcal{C}^{\beta})\quad P^\nu-a.s..\endaligned$$
Moreover, we have $$P^\nu[X\in C([0,\infty),\mathcal{C}^\alpha)]=1.$$
Similar arguments as in the proof of Lemma 3.6 imply that $\forall s, t_0\geq0$
 $$\aligned P^\nu(:p(X(s+t_0)):=\sum_{k=1}^{2N}ka_k\sum_{l=0}^{k-1} C_{k-1}^l\bar{Y}(s,t_0)^l:[Z(s+t_0)+e^{sA}(X(t_0)-Z(t_0))]^{k-1-l}:,\\ X\in C([0,\infty),\mathcal{C}^\alpha),\bar{Y}\in C([0,\infty)^2;\mathcal{C}^{\beta}))=1,\endaligned$$
 where we used that for $s+t_0>0$, $e^{sA}(X(t_0)-Z(t_0))\in \mathcal{C}^\beta$ $P^\nu$-a.s.,  and $:[Z(s+t_0)+e^{sA}(X(t_0)-Z(t_0))]^{k}::=\sum_{l=0}^{k} C_{k}^l(e^{sA}(X(t_0)-Z(t_0)))^l:Z(s+t_0)^{k-l}:\in \mathcal{C}^\alpha$.
 Then using Fubini's theorem we know that
  $$\aligned P^\nu(\int_0^t:p(X(s+t_0)):ds=\sum_{k=1}^{2N}ka_k\sum_{l=0}^{k-1} \int_0^t C_{k-1}^l\bar{Y}(s,t_0)^l:[Z(s+t_0)+e^{sA}(X(t_0)-Z(t_0))]^{k-1-l}:ds,\\   \forall t\geq0, a.e. t_0\geq0, X\in C([0,\infty);\mathcal{C}^\alpha),\bar{Y}\in C([0,\infty)^2;\mathcal{C}^{\beta}))=1.\endaligned$$
  Here we used $X\in C([0,\infty);\mathcal{C}^\alpha)$ for $-\alpha(2N-1)<1$ to make the right hand side of the first equality  meaningful.
 It is obvious that the right hand side of the first equality  is continuous with respect to $t_0$.  Since $\int_0^t:p(X(s+t_0)):ds=\int_{t_0}^{t+t_0}:p(X(s)):ds$ we know that $\int_0^t:p(X(s+t_0)):ds$ is also continuous with respect to $t_0$ and we obtain that $$\aligned P^\nu(\int_0^t:p(X(s+t_0)):ds=\sum_{k=1}^{2N}ka_k\sum_{l=0}^{k-1} \int_0^t C_{k-1}^l\bar{Y}(s,t_0)^l:[Z(s+t_0)+e^{sA}(X(t_0)-Z(t_0))]^{k-1-l}:ds,\\   \forall t,t_0\geq0, X\in C([0,\infty);\mathcal{C}^\alpha),\bar{Y}\in C([0,\infty)^2;\mathcal{C}^{\beta}))=1.\endaligned$$
 This  implies that  there exists a properly $\mathcal{E}$-exceptional set $S_2\supset S_1$ such that for  $z\in \mathcal{C}^{\alpha}\setminus S_2$ under $P^z$
 $$\aligned P^z(X\in C([0,\infty);\mathcal{C}^\alpha), \int_0^t:p(X(s+t_0)):ds=\int_0^t\sum_{k=1}^{2N}ka_k\sum_{l=0}^{k-1} C_{k-1}^l\bar{Y}(s,t_0)^l\\:[Z(s+t_0)+e^{sA}(X(t_0)-Z(t_0))]^{k-1-l}:ds, \forall t, t_0\geq0 )=1.\endaligned$$
 Indeed, define $$\aligned\Omega_0:=&\{\omega:X\in C([0,\infty);\mathcal{C}^\alpha),\int_0^t:p(X(s+t_0)):ds=\int_0^t\sum_{k=1}^{2N}ka_k\sum_{l=0}^{k-1} C_{k-1}^l\bar{Y}(s,t_0)^l\\&:[Z(s+t_0)+e^{sA}(X(t_0)-Z(t_0))]^{k-1-l}:ds, \forall t, t_0\geq0 \},\endaligned$$
 and
 let $\Theta_t:\Omega\rightarrow\Omega, t > 0$, be the canonical shift, i.e. $\Theta_t(\omega) = \omega(\cdot+t),\omega\in \Omega$.
Then it is easy to check that
$$\Theta_t^{-1}\Omega_0\supset \Omega_0,\quad t\in \mathbb{R}^+,$$
and
$$\Omega_0=\bigcap_{t>0, t\in\mathbb{Q}}\Theta_t^{-1}\Omega_0.$$
On the other hand, by the Markov property we know that $$P^z(\Theta_t^{-1}\Omega_0)=P_t(1_{\Omega_0})(z),$$ which by [MR92, Chapter IV Theorem 3.5] is $\mathcal{E}$-quasi-continuous in the sense of [MR92, Chapter III Definition 3.2]  on $E$. It follows that for every $t>0$
$$P^z(\Theta_t^{-1}\Omega_0)=1\quad q.e. z\in E,$$
which yields that
$$P^z(\Omega_0)=1\quad q.e. z\in E.$$
Here q.e. means that there exists a properly $\mathcal{E}$-exceptional set such that outside this exceptional set the result follows. Now $Y$ satisfies (3.4) $P^z$-a.s. for $z\in \mathcal{C}^\alpha\backslash S_2$. Moreover, for $z\in \mathcal{C}^\alpha\backslash S_2$ $Y\in C([0,T];\mathcal{C}^\beta), \bar{Z}\in C([0,T];\mathcal{C}^\alpha) P^z$-a.s., which implies that $$P^z[X(t)\in {\mathcal{C}}^{\alpha}\setminus S_2, \forall t\geq0]=1 \textrm{ for } z\in \mathcal{C}^{\alpha}\setminus S_2.$$
$\hfill\Box$
\vskip.10in

\subsection{Relations between two solutions: starting with solutions to the shifted equation}
Now we fix a stochastic basis $(\Omega,\mathcal{F}, \{\mathcal{F}_t\}_{t\in[0,\infty)},P)$ and on it a cylindrical Wiener process $W$ in $L^2(\mathbb{T}^2)$. Define $\bar{Z}=\int_0^te^{(t-s)A}dW(s)+e^{tA}z$ as in Section 3.2 with $z\in\mathcal{C}^\alpha$ for $\alpha<0$. Now we consider the following equation:
\begin{equation}Y(t)=e^{tA}y-\int_0^te^{(t-s)A}\sum_{k=1}^{2N}ka_k\sum_{l=0}^{k-1} C_{k-1}^lY(s)^l:\bar{Z}(s)^{k-1-l}:ds.\end{equation}
When $N=2$, global existence and uniqueness of the solutions to (3.5) have been obtained in [MW15]. Now we consider general $N\in\mathbb{N}$ and   have the following result. Moreover, by using solutions given by Dirichlet form theory we  also obtain that $\nu$ is an invariant measure of the solution to $\bar{X}=Y_0+\bar{Z}$, where  $Y_0$ is the unique solution to (3.5) with $y=0$.
\vskip.10in
\th{Theorem 3.10} Fix $\alpha,\beta$ such that $0<-\alpha<\beta<\alpha+2$ and  $\beta,-\alpha$ sufficiently small. For  $y\in L^p(\mathbb{T}^2),$ with $p$ even, large enough,  there exists a unique solution to (3.5) in $C((0,T];\mathcal{C}^{\beta})$ equipped with the norm $\sup_{t\in[0,T]}t^{\frac{\beta}{2}+\frac{1}{p}}\|\cdot\|_\beta$.

Moreover, $\nu$ is an invariant measure of the solution to $\bar{X}=Y_0+\bar{Z}$, where  $Y_0$ is the unique solution to (3.5) with $y=0$.

\proof First we  prove  local existence and uniqueness of solutions:
for $y\in L^p(\mathbb{T}^2)$ and a.s. $\omega\in\Omega$ there exists $T^*(z,\omega)$ and a unique solution to (3.5) such that
$$Y(\omega)\in  C((0,T^*(z,\omega)],\mathcal{C}^{\beta}).$$
In fact, we use a fixed point argument in the space $$\mathcal{L}_T:=C((0,T],\mathcal{C}^{\beta})$$
equipped with the norm $\sup_{t\in[0,T]}t^{\frac{\beta}{2}+\frac{1}{p}}\|\cdot\|_\beta$.
By Lemma 3.5 we have  that  $:\bar{Z}^n_\varepsilon:$ converges to $:\bar{Z}^n:$ in  $L^p(\Omega,C((0,T]; {\mathcal{C}}^{\alpha})$ with the norm $\sup_{t\in[0,T]}t^\frac{(\beta-\alpha)n+\rho}{2}\|:\bar{Z}^{n}(t):\|_{\alpha}<\infty$ for $\beta+\alpha>0, \rho>0$. We introduce the following notation $$\|\bar{Z}\|_{\mathfrak{L}}:=\sum_{l=0}^{2N-1}\sup_{\tau\in[0,T]}\tau^{\frac{(\beta-\alpha)l+\rho}{2}}\|:\bar{Z}^{l}(\tau):
\|_{\alpha}.$$
By Lemmas 2.2, 2.3 we obtain that for $Y\in\mathcal{L}_T$
$$\int_0^te^{(t-\tau)A}\sum_{k=1}^{2N}ka_k\sum_{l=0}^{k-1}C_{k-1}^lY^l(\tau):\bar{Z}^{k-1-l}(\tau):d\tau\in\mathcal{L}_T$$
and
\begin{equation}\aligned&\sup_{t\in[0,T]}t^{\beta/2+1/p}\bigg\|\int_0^te^{(t-\tau)A}\sum_{k=1}^{2N}ka_k\sum_{l=0}^{k-1}C_{k-1}^lY^l(\tau):\bar{Z}^{k-1-l}(\tau):d\tau\bigg\|_{\beta}
\\\lesssim&\sup_{t\in[0,T]}t^{\beta/2+1/p}\sum_{k=1}^{2N}\int_0^t(t-\tau)^{\frac{\alpha-\beta}{2}}\sum_{l=0}^{k-1}\|Y(\tau)\|^l_{\beta}
\|:\bar{Z}^{k-1-l}(\tau):\|_{\alpha}
d\tau\\\lesssim&\sup_{t\in[0,T]}t^{\beta/2+1/p}\sum_{k=1}^{2N}\int_0^t(t-\tau)^{\frac{\alpha-\beta}{2}}
(\sum_{l=1}^{k-1}\tau^{-(\frac{\beta}{2}+\frac{1}{p})l-\frac{(\beta-\alpha)(k-1-l)+\rho}{2}}\|Y\|^l_{\mathcal{L}_T}+\tau^{-\frac{(\beta-\alpha)(k-1)+\rho}{2}})\|Z\|_{\mathfrak{L}}
d\tau\\\lesssim&T^{(\alpha-\beta)N+1+\frac{\beta}{2}+\frac{1}{p}-\frac{\rho}{2}}\vee T^{1+\frac{\alpha-\beta-\rho}{2}-(\frac{\beta}{2}+\frac{1}{p})(2N-2)}(\|Y\|^{2N}_{\mathcal{L}_T}+1),\endaligned\end{equation}
with $\beta,-\alpha>0$ small enough and $p>0$ large enough. Here we used Lemma 2.3 in the first inequality and used Lemma 3.5 in the second inequality.
Moreover, by Lemmas 2.1 and 2.2 we have
$$\|e^{tA}y\|_{\beta}\lesssim t^{-(\frac{\beta}{2}+\frac{1}{p})}\|y\|_{L^p(\mathbb{T}^2)}.$$
Similarly we obtain that the iteration mapping is a strict contraction in a bounded ball  $\mathcal{L}_T$ with $T>0$ small enough, which implies the local existence and uniqueness of solutions.
A similar argument as (3.6) also implies that the local solution is  continuous with respect to $(\bar{Z}, :\bar{Z}^2:,...,:\bar{Z}^{2N-1}:)$.

In the following we give an a-priori estimate on the $L^p$ norm of $Y$: Let $Y_\varepsilon$ be the solution to (3.5) with $\bar{Z}$ replaced by $\bar{Z}_\varepsilon$. Then we have
$$\aligned&\frac{1}{p}(\|Y_\varepsilon(t)\|_{L^p}^p-\|Y_\varepsilon(0)\|_{L^p}^p)\\=&\int_0^t[-(p-1)\langle\nabla Y_\varepsilon(s), Y_\varepsilon^{p-2}(s)\nabla Y_\varepsilon(s)\rangle-\|Y_\varepsilon(s)\|_{L^p}^p-\langle \sum_{k=1}^{2N}a_k\sum_{l=0}^{k-1}C_{k-1}^lY_\varepsilon^l(s):\bar{Z}_\varepsilon^{k-1-l}(s):,Y_\varepsilon(s)^{p-1}\rangle ]ds.\endaligned$$
Without loss of generality, suppose that $a_{2N}=\frac{1}{2N}$. Then
$$\aligned&\frac{1}{p}(\|Y_\varepsilon(t)\|_{L^p}^p-\|Y_\varepsilon(0)\|_{L^p}^p)+\int_0^t[(p-1)\langle\nabla Y_\varepsilon(s), Y_\varepsilon(s)^{p-2}\nabla Y_\varepsilon(s)\rangle+\|Y_\varepsilon(s)^{p+2N-2}\|_{L^1}]ds\\&=-\int_0^t[\|Y_\varepsilon(s)\|_{L^p}^p+\langle \Psi(Y_\varepsilon(s),\bar{Z}_\varepsilon(s)),Y_\varepsilon(s)^{p-1}\rangle]ds.\endaligned$$
Here $\Psi(Y_\varepsilon(s),\bar{Z}_\varepsilon(s))=\sum_{k=1}^{2N}ka_k\sum_{l=0}^{k-1}C_{k-1}^lY_\varepsilon^l(s):\bar{Z}_\varepsilon^{k-1-l}(s):
-Y_\varepsilon^{2N-1}(s)$.
Now we only consider $\langle Y_\varepsilon(s)^{2N-2}\bar{Z}_\varepsilon(s), Y_\varepsilon(s)^{p-1} \rangle$. The other terms can be estimated similarly.
We have $$\langle Y_\varepsilon(s)^{2N-2}\bar{Z}_\varepsilon(s), Y_\varepsilon(s)^{p-1} \rangle=\langle Y_\varepsilon(s)^{2N+p-3},\bar{Z}_\varepsilon(s)\rangle.$$
In the following we omit $\varepsilon$ if there's no confusion. Then Lemma 2.3 implies  the following duality
$$|\langle Y(s)^{2N+p-3},\bar{Z}(s)\rangle|\lesssim \|Y(s)^{2N+p-3}\|_{\mathcal{B}^{-\alpha}_{1,1}}\|\bar{Z}(s)\|_{{\alpha}}.$$
Moreover, we have
$$\aligned&\|Y(s)^{2N+p-3}\|_{B^{-\alpha}_{1,1}}\lesssim
\|\Lambda^{\beta_0}Y(s)^{2N+p-3}\|_{L^{p_0}}\lesssim\|\Lambda^{{\beta_0}} Y^{\frac{p}{2}+N-\frac{3}{2}}\|_{L^{p_1}}\|Y^{\frac{p}{2}+N-\frac{3}{2}}\|_{L^{q_1}},\endaligned$$
with $\Lambda=(-\Delta+I)^{1/2}, \beta_0>-\alpha>0, p_0>1,   \frac{1}{p_1}+\frac{1}{q_1}=\frac{1}{p_0}$, where we used Lemma 2.1 in the first inequality and Lemma 2.4 in the second inequality.  Now we estimate each term separately: Lemmas 2.1 and 2.4 imply that
$$\aligned&\|\Lambda^{\beta_0} Y^{\frac{p}{2}+N-\frac{3}{2}}\|_{L^{p_1}}\lesssim\|\Lambda^{\beta_1} Y^{\frac{p}{2}+N-\frac{3}{2}}\|_{L^{p_2}}\lesssim \|\Lambda Y^{\frac{p}{2}+N-\frac{3}{2}}\|_{L^{p_2}}^{\beta_1}\|Y^{\frac{p}{2}+N-\frac{3}{2}}\|_{L^{p_2}}^{1-\beta_1},\endaligned$$
where $\beta_1=\beta_0+\frac{2}{p_2}-\frac{2}{p_1}$, $1<p_2<p_1<2$.
For $\|\Lambda Y^{\frac{p}{2}+N-\frac{3}{2}}\|_{L^{p_2}}$ we have
$$\aligned\|\Lambda Y^{\frac{p}{2}+N-\frac{3}{2}}\|_{L^{p_2}}\lesssim& \| Y^{\frac{p}{2}+N-\frac{5}{2}}\nabla Y\|_{L^{p_2}}\lesssim \| Y^{p-2}|\nabla Y|^2\|_{L^{1}}^{\frac{1}{2}}\|Y^{\frac{p_2(2N-3)}{2-p_2}}\|_{L^1}^{\frac{1}{p_2}-\frac{1}{2}},\endaligned$$
 where we used  H\"{o}lder's inequality in the last inequality. Furthermore, we have
$$\|Y^{\frac{p}{2}+N-\frac{3}{2}}\|_{L^{p_2}}\lesssim\|Y^{\frac{p}{2}+N-\frac{3}{2}}\|_{L^{q_1}}\lesssim \|Y^{p+2N-2}\|_{L^1}^{\frac{\frac{p}{2}+N-\frac{3}{2}}{p+2N-2}},$$
with $2<q_1<\frac{p+2N-2}{{\frac{p}{2}+N-\frac{3}{2}}}$.
Choose $p$ large enough (depending only on $N$) such that
\begin{equation}\frac{4N+p-5}{2(p+2N-2)}<\frac{1}{p_2}<\frac{1}{2(N-1)}+\frac{1}{p_1}<\frac{1}{2(N-1)}+\frac{\frac{p}{2}+N-\frac{1}{2}}{p+2N-2},\end{equation}
which implies that $\frac{p_2(2N-3)}{2-p_2}\leq p+2N-2$. Thus  we obtain that
$$|\langle Y(s)^{2N+p-3},\bar{Z}(s)\rangle|\lesssim \|\bar{Z}(s)\|_{{\alpha}}\|Y^{p+2N-2}\|_{L^1}^{\beta_1\frac{N-\frac{3}{2}}{p+2N-2}+(2-\beta_1)\frac{\frac{p}{2}+N-\frac{3}{2}}{p+2N-2}}\| Y^{p-2}|\nabla Y|^2\|_{L^{1}}^{\frac{\beta_1}{2}}.$$
By (3.7) we  know that for $\beta_0$ small enough, $\beta_1 (N-1)<1$ and we have $$\beta_1\frac{N-\frac{3}{2}}{p+2N-2}+(2-\beta_1)\frac{\frac{p}{2}+N-\frac{3}{2}}{p+2N-2}+\frac{\beta_1}{2}<1,$$
which implies that  there exists $\gamma>1$ such that
$$|\langle Y(s)^{2N+p-3},\bar{Z}(s)\rangle|\lesssim \|\bar{Z}(s)\|_{{\alpha}}^\gamma+(1+\varepsilon\|Y^{p+2N-2}\|_{L^1}+\varepsilon\| Y^{p-2}|\nabla Y|^2\|_{L^{1}}).$$
We can do similar calculations for the other terms in $\Psi(Y(s),\bar{Z}(s))$.  Then we  deduce that there exist $0<\rho_1=\frac{\gamma(2N-1)(\beta-\alpha)+\gamma\rho}{2}<1, \gamma>1$ such that
$$\aligned&|\langle \Psi(Y(s),\bar{Z}(s)), Y(s)^{p-1}\rangle|\\\lesssim& s^{-\rho_1}|\bar{Z}|_{\mathfrak{L}}^\gamma
+\varepsilon (\|Y\|_{L^{p+2N-2}}^{p+2N-2}+\| Y^{p-2}|\nabla Y|^2\|_{L^{1}}).\endaligned$$
Here $\rho_1$ can be chosen less than $1$ since $\beta-\alpha>0$ can be chosen small enough.
Hence  Gronwall's inequality yields a uniform estimate
$$\sup_{t\in[0,T]}\|Y_\varepsilon(t)\|_{L^p}^p\lesssim C_T+\|Y(0)\|_{L^p}^p.$$
Therefore, we  extend the local solution to the unique global solution to (3.5) in $C((0,T], \mathcal{C}^\beta)$.

Moreover, consider $\bar{X}=Y_0+\bar{Z}$, where $Y_0$ is the unique solution to (3.5) with $y=0$. By Theorem 3.9 and the uniqueness of the solution to (3.5) we know that $\bar{X}$ has the same law as the solution $X$ constructed by using Dirichlet forms, which combining with $\nu(\mathcal{C}^\alpha)=1$ implies that $\nu$ is an invariant measure of $\bar{X}$.
$\hfill\Box$
\vskip.10in

In the following we start from the transition semigroup of $\bar{X}$: Let $p_t$ be the transition semigroup (of sub-probability kernels) associated with  $\bar{X}$. Since  $\nu(\mathcal{C}^\alpha)=1$, $p_t$ can be extended to a kernel on $E$ by setting
 $$p_t(z,dy)=\delta_z(dy),$$
 for $z\in E\backslash \mathcal{C}^\alpha$, where $\delta_z$ denotes the Dirac measure in $z$.
 By Theorem 3.10 we have $$\int p_tf\nu(dx)= \int f \nu(dx),$$ for $f\in L^2(E;\nu)$ and $$p_tf\rightarrow^{t\rightarrow0} f,$$ for $f\in C_b(\mathcal{C}^\alpha)$.  By [MR92, Chapter II, Subsection 4a] $(p_t)_{t>0}$ uniquely determines a strongly continuous contraction semigroup $(T_t^1)_{t>0}$ of operators on $L^2(E;\nu)$.
By the pathwise uniqueness of the solution to (3.5) we obtain that  $p_tf(z)=P_tf(z)$ for all $f\in\mathcal{B}_b(E), t>0$ and $z\in \mathcal{C}^\alpha\backslash S_2$, which implies that $p_t$ and hence $T_t^1$ is $\nu$-symmetric. Here $P_t$ is the transition semigroup properly associated with Dirichlet form $(\mathcal{E},D(\mathcal{E}))$ defined in Section 3.2. Then there exists a  Dirichlet form $(\mathcal{E}_1,D(\mathcal{E}_1))$ associated with $T_t^1$. Moreover, $(\mathcal{E}_1,D(\mathcal{E}_1))=(\mathcal{E},D(\mathcal{E}))$, where $(\mathcal{E},D(\mathcal{E}))$ is the Dirichlet form obtained in Section 3.2. Hence for $(\mathcal{E}_1,D(\mathcal{E}_1))$ the results in Theorem 3.8 hold.

\vskip.10in
\subsection{Markov uniqueness in the restricted sense}

In this subsection we prove Markov uniqueness in the restricted sense and the uniqueness of the probabilistically weak solutions to (1.1).

By [MR92, Chap.
4, Sect. 4b] it follows that there is a point separating countable $\mathbb{Q}$-vector
space $D\subset \mathcal{F}C_b^\infty$ such that $D\subset D(L({\mathcal{E}}))$. Let $\mathcal{E}^{\textrm{q.r.}}$  be the set of all quasi-regular Dirichlet forms $(\tilde{\mathcal{E}}, D(\tilde{\mathcal{E}}))$ (cf. [MR92]) on $L^2(E;\nu)$ such that $D\subset D(L(\tilde{\mathcal{E}}))$ and $\tilde{\mathcal{E}}=\mathcal{E}$ on $D\times D$.
Here for a Dirichlet form $(\tilde{\mathcal{E}}, D(\tilde{\mathcal{E}}))$ we denote its generator by $(L(\tilde{\mathcal{E}}), D(L(\tilde{\mathcal{E}})))$.

In the following we  consider the martingale problem in the sense of [AR95] and probabilistically weak solutions to (1.1):
\vskip.10in

\th{Definition 3.11}(i) A $\nu$-special standard process $M=(\Omega,\mathcal{F}, (\mathcal{M}_t), X_t, (P^z))$ in the sense of [MR92, Chapter IV] with state space $E$ is said to solve the martingale problem for $(L(\mathcal{E}),  D)$ if for all $u\in  D$, $u(X(t))-u(X(0))-\int_0^t L(\mathcal{E})u(X(s))ds$, $t\geq0$, is an $(\mathcal{M}_t)$-martingale under $P^\nu$.

 (ii) A $\nu$-special standard process $M=(\Omega,\mathcal{F}, (\mathcal{M}_t), X_t, (P^z))$  with state space $E$ is called a  probabilistically weak solution to (1.1) if there exists a map $W:\Omega\rightarrow C([0,\infty);E)$  such that for $\nu$-a.e. $z$ under $P^z$, $W$ is an $\mathcal{M}_t$- cylindrical Wiener process and the sample paths of the associated  process satisfy (3.3)  for all $l\in H^{2+s}$, $s>0$.
\vskip.10in

\th{Remark} If $M$ is a  probabilistically weak solution to (1.1), we can easily check that it also solves the martingale problem. Conversely, if $M$ solves the martingale problem, by the martingale representation theorem in [O05] there exist a stochastic basis $(\Omega,\mathcal{F}, \{\mathcal{F}_t \}_{t\in[0,\infty)},P)$, a cylindrical Wiener process
$W$ on the space $E$ and a progressively measurable process $X$ such that  $X$ satisfies (3.3)  for $l\in H^{2+s}$, $s>0$.
\vskip.10in
To explain the uniqueness result below we also introduce the following  concept:

Two strong Markov processes $M$ and $M'$ with state space $E$ and transition semigroups $(p_t)_{t>0}$ and $(p_t')_{t>0}$ are called \emph{$\nu$-equivalent} if there exists $S\in\mathcal{B}(E)$ such that (i) $\nu(E\backslash S)=0$, (ii) $P^z[X(t)\in S, \forall t\geq0]=P'^z[X'(t)\in S, \forall t\geq0]=1,  z\in S$, (iii) $p_tf(z)=p_t'f(z)$ for all $f\in\mathcal{B}_b(E), t>0$ and $z\in S$.

\vskip.10in

Combining Theorem 3.9 and Theorem 3.10,  we  obtain Markov uniqueness in the restricted sense for $(L(\mathcal{E}), D)$ (see part (ii)), but in fact prove a much stronger result, i.e. the uniqueness of probabilistically weak solutions to (1.1) (see part (i)):
\vskip.10in

\th{Theorem 3.12} (i) There exists  (up to  $\nu$-equivalence) exactly one  probabilistically weak solution $M$ to (1.1) satisfying  $P^z(X\in C([0,\infty);E))=1$ for $\nu$-a.e.  and has $\nu$ as a subinvariant measure, i.e. for the transition semigroup $(p_t)_{t\geq0}$, $\int p_t fd\nu\leq \int fd\nu$ for $f\in L^2(E;\nu)$.

(ii) $\sharp\mathcal{E}^{\textrm{q.r.}}=1$. Moreover,  there exists (up to $\nu$-equivalence) exactly one $\nu$-special standard process $M$ with state space $E$  associated with a Dirichlet form $(\mathcal{E},D(\mathcal{E}))$ solving the martingale problem for $(L(\mathcal{E}),D)$ .

\proof  For (i), suppose that $M^1$ is a probabilistically weak solution to (1.1) and  let $p_t^1$ be the transition semigroup (of sub-probability kernels) associated with  $M^1$.
Since $\nu$ is a subinvariant measure and $$p_t^1f\rightarrow^{t\rightarrow0} f,$$ for $f\in \mathcal{F}C_b^\infty$, by [MR92, Chapter II, Subsection 4a] $(p_t^1)_{t>0}$ uniquely determines a strongly continuous contraction semigroup $(T_t^1)_{t>0}$ of operators on $L^2(E;\nu)$.
By the proof of Theorem 3.9 we know that the solution to (3.3) having $\nu$ as a subinvariant measure also satisfies (3.4). Moreover, by the pathwise uniqueness of solutions to (3.4) we obtain that  $p_t^1f(z)=P_tf(z)$ $\nu$-a.e. for all $f\in\mathcal{B}_b(E), t>0$, which implies that $p_t^1$ is associated with the Dirichlet form $(\mathcal{E},D(\mathcal{E}))$ obtained in Section 3.2. Here $P_t$ is the semigroup properly associated with $(\mathcal{E},D(\mathcal{E}))$ obtained in Section 3.2. Since $M^1$ is a $\nu$-special standard process and has continuous paths, by [MR92, Chapter 4, Theorem 1.15, Theorem 5.1] $M^1$ is properly associated with $(\mathcal{E},D(\mathcal{E}))$. Then by [MR92, Chapter 4, Theorem 6.4] $M^1$ is $\nu$-equivalent to $M$ obtained in Section 3.2, which implies (i) easily.

 The second result in (ii) follows from the first result and [AR95, Theorem 3.4]. We only prove the first. Since for every $\tilde{\mathcal{E}}\in \mathcal{E}^{\textrm{q.r.}}$ there exists a unique Markov process $\tilde{M}$ associated with $\tilde{\mathcal{E}}$ and  Theorem 3.8 holds for $\tilde{M}$, by Theorems 3.9 and 3.10 we know that for the semigroup $\tilde{p}_t$ associated with $\tilde{M}$ we have $\tilde{p}_tf=P_tf$ $\nu$-a.e. for $f\in \mathcal{B}_b(E)$, which implies that $\tilde{p}_t$ is a $\nu$-version of the semigroup $T_t$ associated with $(\mathcal{E},D(\mathcal{E}))$. Then by [MR92, Chapter I] we know that $(\mathcal{E},D(\mathcal{E}))=(\tilde{\mathcal{E}},D(\tilde{\mathcal{E}}))$. Now (ii) follows.$\hfill\Box$

\subsection{Stationary solution}

Now we consider the stationary case. In this case, we can obtain a probabilistically strong solution to (3.3). Take two different stationary solutions $X_1, X_2$ to (3.3) with the same initial condition $\eta\in \mathcal{C}^{\alpha}$, $\alpha<0$, $-\alpha$ small enough,  having the distribution  $\nu$. We have
$$X_i(t)=e^{tA}\eta-\int_0^te^{(t-\tau)A}:p(X_i(\tau)):d\tau+Z(t),$$
where $Z$ is the stochastic convolution
$$Z(t)=\int_{0}^te^{(t-s)A}d W(s).$$
By a similar argument as  the proof of Theorem 3.9 and using  Lemma 3.1 we have that for every $p>1$
$$E\int_0^T\|:p(X_i(\tau)):\|_{{\alpha}}^pd\tau=T\int\|:p(\phi):\|_{{\alpha}}^p\nu(d\phi)<\infty.$$
Then  Lemma 2.2 implies that for $\alpha<0$,   $-\alpha<\beta<\alpha+2$
$$\int_0^te^{(t-\tau)A}:p(X_i(\tau)):d\tau\in C([0,T];\mathcal{C}^{\beta})\quad P-a.s.. $$
Then by Lemma 2.2 we conclude that
$$X_i-Z\in  C((0,T];\mathcal{C}^{\beta})\quad P-a.s.,$$
where $C((0,T];\mathcal{C}^{\beta})$ is equipped with the norm $\sup_{t\in[0,T]}t^{\frac{\beta-\alpha}{2}}\|\cdot\|_{\beta}$.
Moreover,  similar arguments as in the proof of Theorem 3.9 yield that if $\alpha<0$ with $-\alpha$ small enough, $X_i-Z$ is a solution to the following equation
\begin{equation}Y(t)=e^{tA}\eta+\int_0^te^{(t-s)A}\sum_{k=1}^{2N}ka_k\sum_{l=0}^{k-1} C_{k-1}^lY(s)^l:Z(s)^{k-1-l}:ds.\end{equation}
Here the Wick powers of $Z$ are defined as in Lemma 3.4.

Now by similar calculations as in (3.6) we obtain local uniqueness of the solution to (3.8), which implies that
$$X_1-Z=X_2-Z \textrm{ on } [0,T]\quad P-a.s..$$
Then the pathwise uniqueness holds for the stationary solution to (3.3). Now by the existence of the stationary martingale solution ( cf. [MR99]) and the  Yamada-Watanabe Theorem in [Kur07] we obtain:
 \vskip.10in

\th{Theorem 3.13} For any initial condition $X(0)\in \mathcal{C}^\alpha$ with distribution $\nu$ and $\alpha<0$, $-\alpha$ small enough, there exists a unique probabilistically strong solution $X$ to   (3.3) such that $X$ is a stationary process, i.e.  for every probability space $(\Omega,\mathcal{F},\{\mathcal{F}_t\}_{t\in [0,T]},P)$ with an $\mathcal{F}_t$-Wiener process $W$, there exists  an $\mathcal{F}_t$-adapted process $X:[0,T]\times \Omega\rightarrow E$ such that
for $P-a.e.$ $\omega\in \Omega$ $X$ satisfies (3.3). Moreover, for $0<\beta<\alpha+2$
$$X-Z\in  C((0,T];\mathcal{C}^{\beta})\quad P-a.s..$$

\section{Infinite volume case}

In this section, we analyze the   stochastic quantization equations in infinite volume. The proof is similar as for the finite volume case. However, the  invariant measure $\nu_0$ defined below is more singular and the analysis becomes considerably harder. For simplicity we choose $N=2$. The general case  can be proved similarly. Recall that $\mathcal{S}'(\mathbb{R}^2)$ is the space of tempered Schwartz distributions on $\mathbb{R}^2$ and $\mathcal{S}(\mathbb{R}^2)$ the associated test function space equipped with the usual topology. In this section we fix $\sigma>2$.

\subsection{Wick powers}
Let $\mu_0$ be the mean zero Gaussian measure on $(\mathcal{S}'(\mathbb{R}^2),\mathcal{B}(\mathcal{S}'(\mathbb{R}^2)))$ with covariance
$$\int {}_{\mathcal{S}}\!\langle k_1,z\rangle_{\mathcal{S}'}{}_{\mathcal{S}}\!\langle k_2,z\rangle_{\mathcal{S}'}\mu_0(dz)=\int\int G(x-y)k_1(x)k_2(y)dxdy=:\langle k_1,k_2\rangle_{H_{1}},$$ where $G$ denotes the Green function of the operator $-A$ on $\mathbb{R}^2$.
 \vskip.10in

\textbf{Wick powers on $L^2(\mathcal{S}'(\mathbb{R}^2),\mu_0)$}

Let $H_{1}$ be the real Hilbert space obtained
by completing $S(\mathbb{R}^2)$ w.r.t, the norm associated with the inner product $\langle \cdot,\cdot\rangle_{H_{1}}$. Now for $n\in \mathbb{N}$, let $\mathcal{S}_{-n}$ denote the Hilbert subspace of $\mathcal{S}'(\mathbb{R}^2)$ which is the dual of $\mathcal{S}_n$ defined as the completion of $\mathcal{S}(\mathbb{R}^2)$ w.r.t the norm
$$\|k\|_{\mathcal{S}_n}:=[\sum_{|m|\leq n}\int_{\mathbb{R}^2}(1+|x|^2)^n|(\frac{\partial^{m_1}}{\partial x_1^{m_1}},\frac{\partial^{m_2}}{\partial x_2^{m_2}})k(x)|^2dx]^{1/2}.$$For $h\in H_{1}$ we define $X_h\in L^2(\mathcal{S}'(\mathbb{R}^2),\mu_0)$ by $X_h:=\lim_{n\rightarrow\infty}{ }_{\mathcal{S}}\!\langle k_n,\cdot\rangle_{\mathcal{S}'}\textrm{ in } L^2(\mathcal{S}'(\mathbb{R}^2),\mu_0)$ where $k_n$ is any sequence in $\mathcal{S}(\mathbb{R}^2)$ such that $k_n\rightarrow h$ in $H_{1}$. We have the well-known (Wiener-It\^{o}) chaos decomposition $$L^2(\mathcal{S}'(\mathbb{R}^2),\mu_0)=\bigoplus_{n\geq0}\mathcal{H}_n.$$
For $h\in L^2(\mathbb{R}^2,dx)$ and $n\in \mathbb{N}$, define $:z^n:(h)$ to be the unique element in $\mathcal{H}_n$ such that
$$\int:z^n:(h):\prod_{j=1}^nX_{k_j}:_nd\mu_0=n!\int_{\mathbb{R}^2}\prod_{j=1}^n(\int_{\mathbb{R}^2}G(x-y_j)k_j(y_j)dy_j)h(x)dx$$
where $k_1,...,k_n\in \mathcal{S}(\mathbb{R}^2)$ and $:$ $:_n$  means orthogonal projection onto $\mathcal{H}_n$ (see [S74, V.1] for existence of $:z^n:(h)$).

From now on we  define for $h\in L^2(\mathbb{R}^2,dx)$
$$:P(z):(h):=\frac{1}{4}:z^4:(h).$$We have that $\exp(-:P(z):(h))\in L^p(\mathcal{S}'(\mathbb{R}^2),\mu_0)$ for all $p\in [1,\infty)$ if $h\geq0$ (cf. [AR91, Section 7]), hence the following probability measures (called space-time cut-off quantum fields) are well-defined for $\Lambda\in \mathcal{B}(\mathbb{R}^2), \Lambda$ bounded,
$$\nu_\Lambda:=\frac{\exp{(-:P(z):(1_\Lambda))}}{\int\exp{(-:P(z):(1_\Lambda))}d\mu_0}\mu_0.$$
It has been proven that the weak limit $$\lim_{\Lambda\rightarrow \mathbb{R}^2}\nu_\Lambda=:\nu_0$$ exists as a probability measure on $(\mathcal{S}'(\mathbb{R}^2),\mathcal{B}(\mathcal{S}'(\mathbb{R}^2)))$ having moments of all orders (see [GlJ86] and also [AR91, Section 7]). In particular, it follows by [AR89, Proposition 3.7] that $\nu_0(\mathcal{S}_{-n})=1$ for $n\in \mathbb{N}$ large enough. We emphasize that $\nu_0$ is not absolutely continuous with respect to $\mu_0$ (c.f.[AR91]). By [AR91, Section 7], if $n$ is large enough, there exists a $\mathcal{B}(\mathcal{S}_{-n})/\mathcal{B}(\mathcal{S}_{-n})$-measurable map $:\phi^3::\mathcal{S}_{-n}\rightarrow\mathcal{S}_{-n}$ such that $ _{\mathcal{S}_{-n}}\langle:\phi^3:, l\rangle_{\mathcal{S}_{n}}=:\phi^3:(l)$ $\nu_0$-a.e. for each $l$ with compact support and $\int\|:\phi^3:\|_{\mathcal{S}_{-n}}^2d\nu_0<\infty$.

  For $\phi\in \mathcal{S}'(\mathbb{R}^2)$ define  $$\phi_\varepsilon:=\rho_\varepsilon*\phi$$ with $\rho_\varepsilon$ an approximate delta function,
$$\rho_\varepsilon(x)=\varepsilon^{-2}\rho(\frac{x}{\varepsilon})\in \mathcal{D}, \int \rho=1,$$
and for every $n\in\mathbb{N}$ we set $$:\phi_\varepsilon^n:_C=c_\varepsilon^{n/2}P_n(c_\varepsilon^{-1/2}\phi_\varepsilon),$$
with
$c_\varepsilon=\int\phi^2_\varepsilon\mu_0(d\phi)=\int\int G(x-y)\rho_\varepsilon(y)dy\rho_\varepsilon(x)dx
=\|K_\varepsilon\|_{L^2(\mathbb{R}\times \mathbb{R}^2 )}^2$. Here and in the following   $K(t,x-y)$ is the heat kernel associated with $A$ on $\mathbb{R}^2$ and $K_\varepsilon=K*\rho_\varepsilon$ with $*$ means  convolution in space.  By [GlJ86] we know that for every smooth $g$ with compact support, $\langle:\phi^3_\varepsilon:_C,g\rangle$ converges to $\langle:\phi^3:,g\rangle$ in $L^2(\mathcal{S}'(\mathbb{R}^2),\nu_0)$.

Now we give estimates on the measure $\nu_0$ for  later use.
\vskip.10in

\th{Lemma 4.1} Let $\alpha_0<-\frac{3}{2}$, $\sigma>2$,  $p>1, p\in\mathbb{N}$, then
$$\int\| :\phi^3:\|^{2p}_{\hat{\mathcal{B}}^{\alpha_0,\sigma}_{2p,2p}}\nu_0(d\phi)<\infty.$$

\proof By (2.3) and [GlJ86, Corollary 12.2.4] we have
$$\aligned&\int\| :\phi^3:\|^{2p}_{\hat{\mathcal{B}}^{\alpha_0,\sigma}_{2p,2p}}\nu_0(d\phi)\\\lesssim&\int\bigg(\sum_{n=0}^\infty 2^{2n(\alpha_0-1/p+1)p}\sum_{\psi\in\Psi_\star}\sum_{x\in\Lambda_n}|\langle :\phi^3:,\psi^n_{x}\rangle|^{2p}w(x)\bigg)\nu_0(d\phi)\\\lesssim&\sum_{n=0}^\infty 2^{2n(\alpha_0-1/p+1)p}\sum_{\psi\in\Psi_\star}\sum_{x\in\Lambda_n}\int\langle :\phi^3:,\psi^n_{x}\rangle^{2p}w(x)\nu_0(d\phi)\\\leq&C(p)\sum_{n=0}^\infty 2^{2n(\alpha_0-1/p+1)p}\sum_{\psi\in\Psi_\star}\sum_{x\in\Lambda_n}\|\psi^n_{x}\|_{L^4}^{2p}w(x).\endaligned$$
Recall that the $L^4$-norm of $\psi^n_{x}$ is of order $2^{n/2}$ and that $\Psi$ is a finite set. Thus we obtain that the last term is of order
$$\sum_{n=0}^\infty 2^{2n(\alpha_0+\frac{3}{2})p}\int_{\mathbb{R}^2}w(x)dx.$$
Hence the sums over $n$ and $x$ converge for $\alpha_0<-\frac{3}{2}$. $\hfill\Box$
\vskip.10in
\textbf{Wick powers on a fixed probability space}

 Now we fix a stochastic basis $(\Omega,\mathcal{F},(\mathcal{F}_t)_{t\in[0,\infty)}, P)$ and on it a cylindrical Wiener process $W$ in $L^2(\mathbb{R}^2)$. We have the well-known (Wiener-It\^{o}) chaos decomposition $$L^2(\Omega,\mathcal{F},P)=\bigoplus_{n\geq0}\mathcal{H}'_n.$$
Now for $Z(t)=\int_0^t e^{(t-s)A}dW(s)$,  we can also define Wick powers with respect to different covariances by approximations: Let $Z_\varepsilon(t,y)=\int\int_0^t \langle K_\varepsilon(t-s,y-x),dW(s)\rangle$. Here $\langle\cdot,\cdot\rangle$ means  inner product in $L^2(\mathbb{R}^2)$.
For every $n\in\mathbb{N}$ we set $$:Z_\varepsilon^n(t):_{C_t}=(c_{\varepsilon,t})^{\frac{n}{2}}P_n((c_{\varepsilon,t})^{-\frac{1}{2}}Z_\varepsilon(t))\in \mathcal{H}'_n,$$
where $P_n,n=0,1,...,$ are the Hermite polynomials
and $c_{\varepsilon,t}=\|K_\varepsilon 1_{[0,t]}\|_{L^2(\mathbb{R}\times \mathbb{R}^2)}^2$.

By similar arguments as in the proof of Lemma 3.3 and using (2.3) we have:

 \vskip.10in

 \th{Lemma 4.2} For every $\alpha<0$ and every $p>1$,   $n=2,3$,  $:Z_\varepsilon^n:_{C_t}$ converges in $L^p(\Omega,C([0,T];\hat{\mathcal{B}}^{\alpha,\sigma}_{2p,2p}))$. This limit is called Wick power of $Z(t)$ with respect to the covariance $C_t$ and denoted by $:Z^n(t):_{C_t}$.

 \vskip.10in

By this lemma and a similar argument as in the proof of Lemma 3.4 we can also define $:Z^n(t):_C$.
 \vskip.10in

\th{Lemma 4.3} For every $\alpha<0$ and every $p>1$, $n=2,3, $ $:Z_\varepsilon^n:_{C}=(c_{\varepsilon})^{\frac{n}{2}}P_n((c_{\varepsilon})^{-\frac{1}{2}}Z_\varepsilon)$ converges  in $L^p(\Omega,C((0,T];\hat{\mathcal{B}}^{\alpha,\sigma}_{p,p}))$.  Here the norm for $C((0,T];\hat{\mathcal{B}}^{\alpha,\sigma}_{p,p})$ is  $\sup_{t\in[0,T]}t^\rho\|\cdot\|_{\hat{\mathcal{B}}^{\alpha,\sigma}_{p,p}}$ for $\rho>0$. This limit is called Wick power of $Z(t)$ with respect to covariance $C$ and denoted by $:Z^n(t):_C$. Moreover, for $t>0$
$$:Z^n(t):_C=\sum_{l=0}^{[n/2]}c_t^l\frac{n!}{(n-2l)!l!2^l}:Z^{n-2l}(t):_{C_t},$$
where $c_t:=\lim_{\varepsilon\rightarrow0}(c_{\varepsilon,t}-c_\varepsilon)=-\int_t^\infty\int K(r,x)^2dxdr$, $c_{\varepsilon}=\|K_\varepsilon\|_{L^2(\mathbb{R}\times \mathbb{R}^2)}^2$.

 \vskip.10in

Now we combine the initial value part  with the Wick power by using (3.1). In the following we fix $p_0>3$. For $z\in\hat{\mathcal{B}}^{\alpha,\sigma}_{3p_0,\infty}$ with  $\alpha<0$ we set $V(t)=e^{At}z$, $V_\varepsilon=\rho_\varepsilon* V$,  and
$$\bar{Z}(t)=Z(t)+V(t),\quad \bar{Z}_\varepsilon(t)=Z_\varepsilon(t)+V_\varepsilon(t),$$
$$:\bar{Z}^2(t):_C=:Z(t)^2:_C+V(t)^2+2Z(t)V(t),$$
$$:\bar{Z}^3(t):_C=:Z(t)^3:_C+V(t)^3+3Z(t)V^2(t)+3:Z(t)^2:_CV(t).$$
 $:\bar{Z}^n_\varepsilon:_C$ is defined as $:\bar{Z}^n:_C$ with $Z, V$ replaced by ${Z}_\varepsilon, {V}_\varepsilon$, respectively.
 By Lemma 2.5  we know that $V\in C([0,T];\hat{\mathcal{B}}^{\alpha,\sigma}_{3p_0,\infty})$ and $V\in C((0,T];\hat{\mathcal{B}}^{\beta,\sigma}_{3p_0,\infty})$  for $\beta>\alpha$ with the norm $\sup_{t\in[0,T]}t^{\frac{\beta-\alpha}{2}}\|\cdot\|_{\hat{\mathcal{B}}^{\beta,\sigma}_{3p_0,\infty}}$. Moreover, $$\sup_{t\in[0,T]}t^{\frac{\beta-\alpha}{2}}\|V(t)\|_{\hat{\mathcal{B}}^{\beta,\sigma}_{3p_0,\infty}}\lesssim \|z\|_{\hat{\mathcal{B}}^{\alpha,\sigma}_{3p_0,\infty}}.$$ By Lemma 2.6 we obtain that for  $\alpha<0$, $n=1,2,3,$ $\sigma>2$, $p>1$,  $:\bar{Z}^n:_{C}\in L^p(\Omega,C((0,T];\hat{\mathcal{B}}^{\alpha,\sigma}_{p_0,\infty}))$ and that $:\bar{Z}^n_\varepsilon:_C$ converges to $:\bar{Z}^n:_C$ in  $L^p(\Omega,C((0,T];\hat{\mathcal{B}}^{\alpha,\sigma}_{p_0,\infty}))$. Here the norm for $C((0,T];\hat{\mathcal{B}}^{\alpha,\sigma}_{p_0,\infty})$ is $\sup_{t\in[0,T]}t^\frac{(\beta-\alpha)n+\rho}{2}\|\cdot\|_{\hat{\mathcal{B}}^{\alpha,\sigma}_{p_0,\infty}}$ with $\beta>-\alpha>0, \rho>0$.

\vskip.10in

\textbf{Relations between two different Wick powers}

\vskip.10in
\th{Lemma 4.4} Let $\phi$ be a measurable map from $(\Omega, \mathcal{F},P)$ to $C([0,T],\mathcal{S}^{-n})$ with $n$ large enough, $P\circ \phi(t)^{-1}=\nu_0$ for every $t\in[0,T]$  and let $\bar{Z}$ be defined as above. Assume in addition that $y=\phi-\bar{Z}\in  C([0,T];\hat{\mathcal{B}}^{\beta,\sigma}_{p_0,\infty})$ $P$-a.s. for some $\beta$ with  $\beta+\alpha>0$. Then for every $t>0$
$$:\phi^3(t):=\sum_{k=0}^3C_3^ky^{3-k}(t):\bar{Z}^{k}(t):_{C}\quad P-a.s..$$

\proof By [GlJ86, Theorem 12.2.1] it follows that for every compactly supported smooth function $g$ and $t\geq0$
$$\langle:\phi_\varepsilon(t)^3:_C,g\rangle\rightarrow \langle:\phi(t)^3:,g\rangle \quad \textrm{ in } L^2(\Omega,P).$$
Since $y_\varepsilon=\phi_\varepsilon-\bar{Z}_{\varepsilon}=\rho_\varepsilon*y$, it is obvious that $y_\varepsilon(t)\rightarrow y(t)$ in $\hat{\mathcal{B}}^{\beta-\kappa,\sigma}_{p_0,\infty}$ $P$-a.s. for $\kappa>0, \beta+\alpha-\kappa>0$, which combined with Lemmas 2.6 and 4.3 implies that for $k\in\mathbb{N}$, $k\leq 3$
$$ \langle y_\varepsilon^{3-k}(t):\bar{Z}_{\varepsilon}^{k}:_{C},g\rangle\rightarrow \langle y^{3-k}(t):\bar{Z}^{k}:_{C},g\rangle\quad \textrm{in probability } .$$
Moreover, by (3.1) and similar arguments as the proof of Lemma 3.6 we have
$$\aligned :\phi_\varepsilon(t)^3:_C=\sum_{k=0}^3C_3^k:\bar{Z}_{\varepsilon}^k(t):_C y_\varepsilon^{3-k}(t),
\endaligned$$
which implies $$\langle:\phi(t)^3:,g\rangle=\sum_{k=0}^3C_3^k\langle y^{3-k}(t):\bar{Z}^{k}(t):_{C},g\rangle\quad P-a.s..$$
 by letting $\varepsilon\rightarrow0$. Now the results follow because the test function space is separable.
$\hfill\Box$

In the following, we only use Wick powers $:\cdot:_C$ and we write  $:\cdot:$ for simplicity.

 \vskip.10in
\subsection{Relations between the two solutions}
\vskip.10in

{\textbf{Solutions given by Dirichlet forms}}
\vskip.10in

Now choose $H=L^2(\mathbb{R}^2)$ and $E=\mathcal{S}_{-n}$ for some $n$ large enough. We  define the  Dirichlet form as in [AR91].
Define $$\mathcal{F}C_b^\infty=\{u:u(z)=f({ }_{E^*}\!\langle l_1,z\rangle_E,{ }_{E^*}\!\langle l_2,z\rangle_E,...,{ }_{E^*}\!\langle l_m,z\rangle_E),z\in E, l_1,l_2,...,l_m\in E^*, m\in \mathbb{N}, f\in C_b^\infty(\mathbb{R}^m)\},$$
where $E^*$ denotes the dual space of $E$.
 Define for $u\in \mathcal{F}C_b^\infty$ and $l\in H$, $$\frac{\partial u}{\partial l}(z):=\frac{d}{ds}u(z+sl)|_{s=0},z\in E,$$
  that is, by the chain rule,
  $$\frac{\partial u}{\partial l}(z)=\sum_{j=1}^m\partial_jf({ }_{E^*}\!\langle l_1,z\rangle_E,{ }_{E^*}\!\langle l_2,z\rangle_E,...,{ }_{E^*}\!\langle l_m,z\rangle_E)\langle l_j,l\rangle.$$
Let $Du$ denote the $H$-derivative of $u\in \mathcal{F}C_b^\infty$, i.e. the map from $E$ to $H$ such that $$\langle Du(z),l\rangle=\frac{\partial u}{\partial l}(z)\textrm{ for all } l\in H, z\in E.$$
By [AR91] we easily deduce that the form
$$\mathcal{E}(u,v):=\frac{1}{2}\int_E \langle Du, Dv\rangle_{H}d\nu_0; u,v\in\mathcal{F}C_b^\infty$$
is closable and its closure $(\mathcal{E},D(\mathcal{E}))$ is a quasi-regular Dirichlet form on $L^2(E;\nu_0)$ in the sense of [MR92].  By [AR91, Theorem 3.6] we know that there exists a (Markov) diffusion process $M=(\Omega,\mathcal{F},(X(t))_{t\geq0},(P^z)_{z\in E})$ on $E$ {properly associated with} $(\mathcal{E},D(\mathcal{E}))$.

By [AR91, Theorem 7.11] we have the following:
 \vskip.10in

\th{Theorem 4.5} For each $l$ smooth and compactly supported, we have for the partial log derivative of $\nu_0$ $$\beta_l(z)= -\langle:z^{3}:, l\rangle+_{{\mathcal{S}_{n}}}\!\langle \Delta l-l,z\rangle_{{\mathcal{S}_{-n}}}.$$
If $n$ is large enough there exists a $\mathcal{B}(\mathcal{S}_{-n})/\mathcal{B}(\mathcal{S}_{-n})$-measurable map $\beta:\mathcal{S}_{-n}\rightarrow\mathcal{S}_{-n}$ such that $ _{\mathcal{S}_{-n}}\langle\beta, l\rangle_{\mathcal{S}_{n}}=\beta_l$ $\nu_0$-a.e. for each $l$ with compact support and $\int\|\beta\|_{\mathcal{S}_{-n}}^2d\nu_0<\infty$.

 \vskip.10in
Moreover, by [AR91, Theorem 6.1] we obtain the following results:

\vskip.10in
\th{Theorem 4.6} There exist a map $W:\Omega\rightarrow C([0,\infty);E)$ and a \emph{properly  $\mathcal{E}$-exceptional set} $S\subset E$, i.e. $\nu_0(S)=0$ and $P^z[X(t)\in E\setminus S, \forall t\geq0]=1$ for $z\in E\backslash S$, such that $\forall z\in E\backslash S$ under $P^z$, $W$ is an $\mathcal{M}_t$- cylindrical Wiener process and the sample paths of the associated  process $M=(\Omega,\mathcal{F},(X(t))_{t\geq0},(P^z)_{z\in E})$ on $E$ satisfy the following:  for $l\in {\mathcal{S}_{n}}$ with compact support
 \begin{equation}\aligned{ }_{E^*}\!\langle l,X(t)-X(0)\rangle_{E}=&\int_0^t\langle l,dW(r)\rangle+\int_0^t(-{ }_{{\mathcal{S}_{-n}}}\!\langle :X(r)^{3}:,l\rangle_{{\mathcal{S}_{n}}}\\&+_{{\mathcal{S}_{n}}}\!\langle \Delta l-l,X(r)\rangle_{{\mathcal{S}_{-n}}} )dr\quad \forall t\geq 0 \quad  \textrm{ }P^z\rm{-a.s.}.\endaligned\end{equation}
 Moreover, $\nu_0$ is an invariant measure for $X$ in the sense that $\int p_tu d\nu=\int ud\nu$ for $u\in L^2(E;\nu)\cap\mathcal{B}_b(E)$, where $p_t$ is the transition semigroup for $M$.
 \vskip.10in

\textbf{Relations between the two solutions}

 In the following we discuss the relations between $M$ constructed above and the shifted equation.  For $W$ constructed in Theorem 4.6, define
 $\bar{Z}(t):= \int_{0}^te^{(t-s)A}d W(s)+e^{tA}X(0)$. First we prove the following property for $\nu_0$ by using Theorem 4.6. .
 \vskip.10in
\th{Theorem 4.7} For every $\alpha<0,p\geq1$ we have $$\nu_0(\hat{\mathcal{B}}^{\alpha,\sigma}_{p,\infty})=1.$$
\proof By using [GlJ86, Corollary 12.2.4], a similar argument as in the proof of Lemma 4.1 implies that for $\alpha_1<-1/2,p>1$
$$\nu_0(\hat{\mathcal{B}}^{\alpha_1,\sigma}_{2p,2p})=1.$$
Then by Theorem 4.6 we have that for $z\in \hat{\mathcal{B}}^{\alpha_1,\sigma}_{p,\infty}\cap (E\setminus S)$ with $p>1$ under $P^z$

$$X(t)=-\int_0^te^{(t-\tau)A}:X(\tau)^3:d\tau+\bar{Z}(t).$$
By a similar calculation as in the proof of Lemma 3.3, it follows that for every $\alpha<0, p>1, t>0,$ $$\int_{0}^te^{(t-s)A}d W(s)\in \hat{\mathcal{B}}^{\alpha,\sigma}_{p,\infty}\quad P^{\nu_0}-a.s..$$
Moreover,  by Lemma 2.5 for $t>0$ $$\|e^{tA}z\|_{\hat{\mathcal{B}}^{\alpha,\sigma}_{p,\infty}}\lesssim t^{-\frac{(\alpha-\alpha_1)}{2}}\|z\|_{\hat{\mathcal{B}}^{\alpha_1,\sigma}_{p,\infty}},$$
which implies that for every $t>0,\alpha<0$, $$\bar{Z}(t)\in \hat{\mathcal{B}}^{\alpha,\sigma}_{p,\infty}\quad P^{\nu_0}-a.s..$$
Since $\nu_0$ is an invariant measure for $M$, by Lemma 4.1 we conclude that  for every $\alpha_0<-\frac{3}{2},T>0$, $p>1$,
$$ E^{\nu_0}\int_0^T\|:X(\tau)^3:\|_{\hat{\mathcal{B}}^{\alpha_0,\sigma}_{p,p}}^pd\tau=T\int\|:\phi^3:\|_{\hat{\mathcal{B}}^{\alpha_0,\sigma}_{p,p}}^p\nu_0(d\phi)<\infty.$$
 Then by  Lemma 2.5 we have that for $p>\frac{2}{2-(\alpha-\alpha_0)}$
 $$\aligned &E^{\nu_0}\sup_{t\in[0,T]}\|\int_0^te^{(t-\tau)A}:X(\tau)^3:d\tau\|_{\hat{\mathcal{B}}^{\alpha,\sigma}_{p,p}}\lesssim E^{\nu_0}\sup_{t\in[0,T]}\int_0^t(t-\tau)^{-\frac{\alpha-\alpha_0}{2}}\|:X(\tau)^3:\|_{\hat{\mathcal{B}}^{\alpha_0,\sigma}_{p,p}}d\tau
 \\\lesssim &\bigg(E^{\nu_0}\int_0^T\|:X(\tau)^3:\|_{\hat{\mathcal{B}}^{\alpha_0,\sigma}_{p,p}}^pd\tau\bigg)^{\frac{1}{p}}<\infty\endaligned.$$
 Here in the last inequality we used H\"{o}lder's inequality. Thus, by Lemma 2.6 for every $t>0$ $X(t)\in \hat{\mathcal{B}}^{\alpha,\sigma}_{p,\infty}$ $P^{\nu_0}$-a.s., which implies the result since $P^{\nu_0}\circ X(t)^{-1}=\nu_0$.  $\hfill\Box$
 \vskip.10in

 Now we prove that $X-\bar{Z}$ satisfies the shifted equation.

\vskip.10in
 \th{Theorem 4.8} Let $\alpha\in (-\frac{1}{3},0)$ and $p_0>3$. There exists a properly $\mathcal{E}$-exceptional set $S_2\supset S$ in the sense of Theorem 4.6 such that for  $z\in \hat{\mathcal{B}}^{\alpha,\sigma}_{3p_0,\infty}\cap(E\setminus S_2)$   $Y:=X-\bar{Z}\in C([0,T];\hat{\mathcal{B}}^{\beta,\sigma}_{p,p})$ $P^z$-a.s. for every $\beta\in (0,\frac{1}{2}), p>1$,  is a solution to the following equation:
\begin{equation}Y(t)=-\int_0^te^{(t-s)A}\sum_{l=0}^3 C_3^lY(s)^l:\bar{Z}^{3-l}(s):ds.\end{equation}
Moreover,
\begin{equation}P^z[X(t)\in \hat{\mathcal{B}}^{\alpha,\sigma}_{3p_0,\infty}\cap(E\setminus S_2), \forall t\geq0]=1 \textrm{ for } z\in \hat{\mathcal{B}}^{\alpha,\sigma}_{3p_0,\infty}\cap(E\setminus S_2).\end{equation}
 \proof  By Theorem 4.6
 we have that for $z\in E\setminus S$
$$X(t)=-\int_0^te^{(t-\tau)A}:X(\tau)^3:d\tau+\bar{Z}(t)\quad P^z-a.s..$$
Since $\nu_0$ is an invariant measure for $X$, by Lemma 4.1 we conclude that  for every $\alpha_0<-\frac{3}{2},$ $p>1$,
$$\int E^z\int_0^T\|:X(\tau)^3:\|_{\hat{\mathcal{B}}^{\alpha_0,\sigma}_{p,p}}^pd\tau\nu_0(dz)=T\int\|:\phi^3:\|_{\hat{\mathcal{B}}^{\alpha_0,\sigma}_{p,p}}^p\nu_0(d\phi)<\infty,$$
which implies that  there exists a properly $\mathcal{E}$-exceptional set $S_1\supset S$ such that for  $z\in E\setminus S_1$
$$ E^z\int_0^T\|:X(\tau)^3:\|_{\hat{\mathcal{B}}^{\alpha_0,\sigma}_{p,p}}^pd\tau<\infty.$$
By Lemma 2.5 we know that for $0<\beta<\alpha_0+2$  and for $z\in E\backslash S_1$, $p>1$,
$$\int_0^te^{(t-\tau)A}:X(\tau)^3:d\tau\in C([0,T];{\hat{\mathcal{B}}^{\beta,\sigma}_{p,p}})\quad P^z- a.s.. $$
Then we conclude that for every  $z\in E\setminus S_1$, $p>1$
$$X-\bar{Z}\in  C([0,T];{\hat{\mathcal{B}}^{\beta,\sigma}_{p,p}})\quad P^z- a.s..$$
By Theorem 4.7 and the fact that $Z\in C([0,T];\hat{\mathcal{B}}^{\alpha,\sigma}_{3p_0,\infty})$ $P^{\nu_0}$-a.s. for $\alpha\in (-\frac{1}{3},0)$, we obtain that
$$\bar{Z}\in C([0,T];\hat{\mathcal{B}}^{\alpha,\sigma}_{3p_0,\infty})\quad P^{\nu_0}-a.s.. $$
Thus, Lemma 4.4 and similar arguments as in the proof of Theorem 3.9 imply  that
 $$\aligned &P^{\nu_0}[X\in C([0,\infty),\hat{\mathcal{B}}^{\alpha,\sigma}_{3p_0,\infty}), X-\bar{Z}\in C([0,\infty),\hat{\mathcal{B}}^{\beta,\sigma}_{p,p}),  \\&\int_0^t:X(s)^3:ds=\int_0^t\sum_{l=0}^{3} C_3^l(X(s)-\bar{Z}(s))^l:\bar{Z}(s)^{3-l}:ds, \forall t\geq0 ]=1,\endaligned$$
 which combined with similar arguments as in the proof of Theorem 3.9 implies that  there exists a properly $\mathcal{E}$-exceptional set $S_2\supset S$ such that for  $z\in \hat{\mathcal{B}}^{\alpha,\sigma}_{3p_0,\infty}\cap(E\setminus S_2)$
 $$\aligned &P^z[X\in C([0,\infty),\hat{\mathcal{B}}^{\alpha,\sigma}_{3p_0,\infty}), X-\bar{Z}\in C([0,\infty),\hat{\mathcal{B}}^{\beta,\sigma}_{p,p}),  \\&\int_0^t:X(s)^3:ds=\int_0^t\sum_{l=0}^{3} C_3^l(X(s)-\bar{Z}(s))^l:\bar{Z}(s)^{3-l}:ds, \forall t\geq0 ]=1.\endaligned$$
 Now we can conclude the first result. (4.3) follows from the fact that $\bar{Z}\in C([0,T];{\hat{\mathcal{B}}^{\alpha,\sigma}_{3p_0,\infty}})$ and $X-\bar{Z}\in C([0,T];{\hat{\mathcal{B}}^{\beta,\sigma}_{p,p}})$ for every $p>1$.
 $\hfill\Box$

\vskip.10in

Now we deduce the uniqueness of the solution to (4.2) and that $\nu_0$ is  an invariant measure of the solution  $\bar{X}=Y_0+\bar{Z}$, where  $Y_0$ is the unique solution to (4.2).

\vskip.10in
\th{Theorem 4.9} For $0<\beta<\frac{1}{2}$, $\sigma>2$   there exists a unique solution to (4.2) in $C([0,T];\hat{\mathcal{B}}^{\beta,\sigma}_{p,\infty})$.

Moreover, $\nu_0$ is an invariant measure of the solution  $\bar{X}=Y_0+\bar{Z}$, where  $Y_0$ is the unique solution of (4.2).

\proof The first result follows from [MW15, Theorem 9.5] and the second follows from Theorem 4.7 and similar arguments as in  the proof of Theorem 3.10. $\hfill\Box$
\vskip.10in

Similarly as in Section 3.3 we start from the transition semigroup of $\bar{X}$ and  can prove that the  Dirichlet form associated with this transition semigroup is $(\mathcal{E},D(\mathcal{E}))$  obtained before.

\vskip.10in
\subsection{Markov uniqueness in the restricted sense}

All the definitions introduced in Section 3.4 can be transferred here. Combining Theorem 4.8 and Theorem 4.9, we  obtain Markov uniqueness in the restricted sense for $(L(\mathcal{E}), D)$ in the infinite volume case and the uniqueness of probabilistically weak solutions to (1.1):
\vskip.10in

\th{Theorem 4.10}
(i) There exists  (up to  $\nu_0$-equivalence) exactly one  $\nu_0$-special standard process $M$ with state space $E$ which satisfies (4.1) $P^z$-a.s. and $P^z(X\in C([0,\infty);E))=1$ for $\nu_0$-a.e. $z\in E$ and has $\nu_0$ as a subinvariant measure, i.e. for the transition semigroup $(p_t)_{t\geq0}$, $\int p_t fd\nu_0\leq \int fd\nu_0$ for $f\in L^2(E;\nu_0)$.

(ii) $\sharp\mathcal{E}^{\textrm{q.r.}}=1$. Moreover,  there exists (up to $\nu_0$-equivalence) exactly one $\nu_0$-special standard process $M$ with state space $E$ which is associated with a Dirichlet form $(\mathcal{E},D(\mathcal{E}))$ solving the martingale problem for $(L(\mathcal{E}),D)$ .

\proof It follows essentially from the same argument as the proof of Theorem 3.12 and (4.3). $\hfill\Box$


\begin{thebibliography}{99}
\bibitem[AKR12]{}S. Albeverio, H. Kawabi, M. R\"{o}ckner, Strong uniqueness for both Dirichlet operators and stochastic dynamics to Gibbs measures on a path space with exponential interactions, Journal of Functional Analysis 262, 2, 15, 2012,  602-638
    \bibitem[AKR97]{}S. Albeverio, Y. G. Kondratiev, M. R\"{o}ckner, Ergodicity for the Stochastic Dynamics of Quasi-invariant Measures with Applications to Gibbs States, Journal of Functional Analysis, 1997, 149(2), 415-469
\bibitem[AR89]{}S. Albeverio, M. R\"{o}ckner,  Classical Dirichlet forms on topological vector spaces
- Construction of an associated diffusion process. Probab. Th. Ret. Fields. 83,
405-434 (1989)
\bibitem[AR91]{}S. Albeverio, M. R\"{o}ckner, Stochastic differential equations in infinite
dimensions: Solutions via Dirichlet forms, Probab. Theory Related Field 89 (1991) 347-386.
 \bibitem[AR95]{}S. Albeverio, M. R\"{o}ckner,  Dirichlet form methods for uniqueness of martingale problems and applications. Stochastic analysis (Ithaca, NY, 1993), 513-528, Proc. Sympos. Pure Math., 57, Amer. Math. Soc., Providence, RI, 1995
\bibitem[ARZ93a]{}S. Albeverio, M. R\"{o}ckner, T.S. Zhang: Markov uniqueness and its applications to martingale
problems, stochastic differential equations and stochastic quantization . C.R. Math. Rep. Acad.
Sci. Canada XV, 1-6 (1993).
\bibitem[ARZ93b]{} S. Albeverio, M. R\"{o}ckner, T.S. Zhang: Markov uniqueness for a class of infinite dimensional
Dirichlet operators. In: Stochastic Processes and Optimal Control. Stochastic Monographs 7
(eds. H.J.Engelbert et al.) 1-26, Gordon Breach, 1993
\bibitem[BCD11]{} H. Bahouri, J.-Y. Chemin, R. Danchin,  Fourier analysis and nonlinear
partial differential equations, vol. 343 of Grundlehren der Mathematischen
Wissenschaften [Fundamental Principles of Mathematical Sciences]. Springer, Heidelberg,
2011.
\bibitem[DD03]{}G. Da Prato, A. Debussche, Strong solutions to the stochastic quantization equations, Ann.
Probab., 31(4):1900-1916, (2003)
\bibitem[D04]{}G. Da Prato. Kolmogorov Equations for Stochastic PDEs. Birkhäuser, Basel, 2004
\bibitem[GIP13]{} M. Gubinelli, P. Imkeller, N. Perkowski, Paracontrolled distributions and singular PDEs, arXiv:1210.2684
\bibitem[GlJ86]{}  J.  Glimm, A. Jaffe : Quantum physics: a functional integral point of view. New
York Heidelberg Berlin: Springer (1986)
\bibitem[GRS75]{}F. Guerra, J. Rosen, B. Simon: The $P(\Phi)_2$ Euclidean quantum field theory
as classical statistical mechanics. Ann. Math. 101, l 11-259 (1975)
\bibitem[Hai14]{} M. Hairer, A theory of regularity structures. Invent. Math. (2014).
\bibitem[JLM85]{} G. Jona-Lasinio and P. K. Mitter. On the stochastic quantization of field theory. Comm.
Math. Phys., 101(3):409-436, 1985.
\bibitem[Kur07]{}T. G. Kurtz, The Yamada-Watanabe-Engelbert theorem for general stochastic equations and inequalities, \emph{Electronic Journal of Probability}. \textbf{12} (2007), 951-965
    \bibitem[KR07]{}H. Kawabi and M. R\"{o}ckner: Essential self-adjointness of Dirichlet operators on a
path space with Gibbs measures via an SPDE approach, J. Funct. Anal. 242 (2007),
486-518.
\bibitem[LR98]{}V. Liskevich and M. R¨ockner: Strong uniqueness for certain infinite-dimensional
Dirichlet operators and applications to stochastic quantization, Ann. Scuola Norm.
Sup. Pisa Cl. Sci., Serie IV, 27 (1998), no. 1, 69-91.
    \bibitem[MR99]{}R. Mikulevicius, B. Rozovskii, Martingale problems for stochasic PDE's. In Stochastic
partial differential equations: six perspectives, volume 64 of Math. Surveys Monogr.
243-325. Amer. Math. Soc., Providence, RI, 1999
\bibitem[MR92] { }Z. M. Ma, and M. R\"{o}ckner, "Introduction to the theory of (non-symmetric) Dirichlet forms," Springer-Verlag, Berlin/Heidelberg/New York, 1992
 \bibitem[MW15]{}J. Mourrat, H. Weber,  Global well-posedness of the dynamic   $\Phi^4$ model in the plane, arXiv:1501.06191v1
 \bibitem[O05]{} M. Ondrej\'{a}t, Brownian representations of cylindrical local martingales, martingale problem and strong markov property of weak solutions of spdes in Banach spaces, \emph{Czechoslovak Mathematical Journal} \textbf{55} (130)(2005), 1003-1039
  \bibitem[PW81]{} G. Parisi,  Y. S. Wu. Perturbation theory without gauge fixing. Sci. Sinica 24,
no. 4, (1981), 483–496.
 \bibitem[Re95]{} S. Resnick, Danamical Problems in Non-linear Advective Partial Differential Equations, PhD thesis, University of Chicago, Chicago (1995)
    \bibitem[R86]{ }M. R\"{o}ckner, Specifications and Martin boundaries for $P(\phi)_2$-random fields.
\emph{ Commun. Math. Phys.}\textbf{106}, 105-135 (1986)
\bibitem[RZ94] { }M. R\"{o}ckner and T.S. Zhang, Uniqueness of Generalized Schr\"{o}dinger Operators
and Applications, Part II \emph{Journal of Functional Analysis}. \textbf{119} (1994), 455-467
\bibitem[RZZ12] { }M. R\"{o}ckner, R. Zhu, X. Zhu, The stochastic reflection problem on an infinite dimensional convex set and BV functions in a Gelfand triple,
The Annals of Probability
, Vol. 40, No. 4, 1759-1794, (2012)
\bibitem[RZZ15] { }M. R\"{o}ckner, R. Zhu, X. Zhu, BV functions in a Gelfand triple for differentiable measure and its applications
,  Forum Mathematicum.  27, 3, 1657-1687(2015)
\bibitem[RZZ15a] { }M. R\"{o}ckner, R. Zhu, X. Zhu, Sub and supercritical stochastic quasi-geostrophic equation, The Annals of Probability
2015,  43,  3, 1202-1273
\bibitem[S74]{}B. Simon, The $P(\phi)_2$ Euclidean (Quantum) field theory. Princeton: Princeton
University Press (1974)
\bibitem[S85]{}W. Sickel, Periodic spaces and relations to strong summability of multiple Fourier series. Math. Nachr.
124, 15-44 (1985)

\bibitem[SW71]{}E. M. Stein, G. L. Weiss, Introduction to Fourier Analysis on Euclidean Spaces,
Princeton University Press, 1971
\bibitem[Tri78]{} H. Triebel,  Interpolation theory, function spaces, differential operators. North-Holland Mathematical
Library 18. North-Holland Publishing Co. Amsterdam-New York 1978.
\bibitem[Tri83]{}H. Triebel, Theory of function spaces. Basel, Birkh\"{a}user, (1983)
\bibitem[Tri06]{}H. Triebel, Theory of function spaces III. Basel, Birkh\"{a}user, (2006)
\bibitem[ZZ15] { }R. Zhu, X. Zhu, A Wong-Zakai theorem for $\phi^4_3$ model, arXiv:1504.04143



\end{thebibliography}
\end{document}